\documentclass[11pt]{article}
\usepackage{graphicx} % Required for inserting images

\usepackage{mathrsfs}
\usepackage[margin=1in]{geometry}
\usepackage{amsbsy,amsmath,amsthm,amssymb,graphicx,subcaption}
\usepackage{algorithm}
\usepackage{algpseudocode}
\usepackage{verbatim}
\usepackage{natbib}
\usepackage{microtype}
\usepackage{xcolor}
\usepackage{hyperref}

\usepackage{makecell,booktabs}

\makeatletter
\def\namedlabel#1#2{\begingroup
	#2%
	\def\@currentlabel{#2}%
	\phantomsection\label{#1}\endgroup
}
\makeatother

\newcommand{\df}{\mathrm{d}}
\newcommand{\X}{\mathsf{X}}
\newcommand{\Y}{\mathsf{Y}}
\newcommand{\Z}{\mathsf{Z}}

\newcommand{\ind}{\mathbf{1}}

\newtheorem{theorem}{Theorem}
\newtheorem{lemma}[theorem]{Lemma}
\newtheorem{corollary}[theorem]{Corollary}
\newtheorem{proposition}[theorem]{Proposition}

\newtheorem{remark}[theorem]{Remark}

\title{{\bf Decomposition Bounds for MCMC with Mixture Targets: Multimodality and Heavy-Tailedness}}
\author{Qian Qin \\ School of Statistics \\ University of Minnesota}
\date{}

\begin{document}
	
	\maketitle
	
	\begin{abstract}
		We establish quantitative mixing-time bounds for Markov chain Monte Carlo algorithms targeting multimodal and heavy-tailed mixtures. 
		The dynamics of a Markov chain are decomposed into local and global levels, and possibly weak Poincar\'{e} inequalities at both levels are combined to form convergence bounds.
		We apply this framework to location and scale mixtures of strongly log-concave distributions, obtaining explicit convergence bounds for random-walk Metropolis, proximal, and random-scan Gibbs samplers. 
		For location mixtures connected along a finite path, the bounds depend polynomially on the number of components and exponentially on neighboring separation, remaining informative as the overall diameter grows. 
		For scale mixtures, we relate the global dynamics of the chain to the isoperimetry of the log-scale mixing distribution, and illustrate the effect of heavy-tailed targets by examining inverse-gamma mixing distributions.
	\end{abstract}
	
	\section{Introduction} \label{sec:intro}

	Mixture representations provide a way to separate two sources of slow convergence in Markov chain Monte Carlo (MCMC): slow exploration within individual components and limited communication between components. 
	Understanding their combined effect is useful for targets that are multimodal or heavy-tailed. 
	We develop quantitative convergence bounds that distinguish these mechanisms for random-walk Metropolis and data augmentation algorithms, including proximal and random-scan Gibbs samplers. 
	Our applications concern location and scale mixtures of smooth, strongly log-concave distributions.
	
	Mixture decompositions separate the dynamic of a Markov chain into local and global components.
	They have been used extensively for annealing and tempering methods \citep{atchade2021approximate,lee2018beyond,ge2018simulated,guo2024provable,lee2024convergence,zhou2025polynomial,garg2025restricted}. 
	\citet{ge2018simulated} establish polynomial-time sampling guarantees for simulated tempering targeting mixtures of translates of a common strongly log-concave density. 
	Subsequent work further extends the framework \citep{garg2025restricted,zhou2025polynomial}. These results demonstrate that temperature transitions can substantially improve convergence when components are widely separated. In particular, \citet[Corollary~4.3 and Remark~4.4]{zhou2025polynomial} obtains a bound with polylogarithmic dependence on the maximal component separation for simulated tempering.

	We use this approach to obtain explicit mixing-time bounds for the non-tempering-based algorithms mentioned at the beginning of the article.
	We provide general bounds for location mixtures and scale mixtures of strongly log-concave distributions.
	Our main methodological result combines possibly weak Poincar\'{e} inequalities for component kernels with a possibly weak Poincar\'{e} inequality for an auxiliary kernel on the mixture-index space. 
	This extends the soft-decomposition approach and related weak Poincar\'{e} decomposition results \citep{ge2018simulated,huang2025weak}, and provides a common framework for geometric and subgeometric convergence bounds.
	
	For location mixtures, the mixing-time bounds quantify the connectivity among components based on their relative locations. 
	We study an illustrative case that consists of $p$ equally weighted components whose locations form a path, with neighboring distances bounded by $\kappa$. The resulting bounds contain a geometric factor $p^3\exp(L\kappa^2/4)$, where $L$ bounds the curvature of the component potential. Thus, the overall diameter may grow while the bounds remain polynomial in the number of components, provided neighboring separation remains controlled.
	This example identifies a regime in which chains of overlapping components support overall convergence across a large region.
	
	For scale mixtures, the isoperimetry of the log-scale mixing distribution is used to quantify the connectivity among the components.
	An inverse-gamma mixing distribution illustrates how the bounds deteriorate as the tails become increasingly heavy.
	
	To set the stage of our study, we now define the MCMC algorithms and mixture targets that will be focused on.

	Let $f_X$ be a probability distribution on a Polish space~$\X$ with Borel $\sigma$-algebra~$\mathcal{B}$. Let $K:\X\times\mathcal{B}\to[0,1]$ be a Markov transition kernel (Mtk) reversible with respect to $f_X$.
	The kernel~$K$ is assumed to correspond to some MCMC algorithm.
	In particular, the following examples are studied.
	
	\begin{itemize}
		\item \namedlabel{ex:MH}{\textsf{Alg-RWMH}}:
		Let $\X = \mathbb{R}^d$.
		Suppose that $f_X$ has density $Df_X$, and let $q(x,\cdot), \, x \in \mathbb{R}^d$ be the density of the $\mathrm{N}_d(x,\sigma^2I_d)$ distribution for some $\sigma > 0$.
		The Markov transition kernel of a random-walk Metropolis-Hastings algorithm is
		\[
		K(x,\mathsf{A})=\int_{\mathsf{A}} \operatorname{acc}(x,x') \, q(x,x')\,\df x'
		+\left[1-\int_{\X} \operatorname{acc}(x,x') \, q(x,x')\,\df x'\right]\ind_{\mathsf{A}}(x),
		\]
		where
		\[
		\operatorname{acc}(x,x')=1\wedge\frac{Df_X(x')q(x',x)}{Df_X(x)q(x,x')}.
		\]
		Here, $\wedge$ denotes the minimum. Given the current state $x$, the algorithm proposes $x'$ from $q(x,\cdot)$ and accepts it with probability $\operatorname{acc}(x,x')$; otherwise it retains $x$.
		
		\item Data augmentation (DA) algorithm:
		Let $\Y$ be a Polish space with its Borel $\sigma$-algebra. Given a joint distribution
		$f_{X,Y}(\df x,\df y)=f_X(\df x)f_{Y\mid X}(\df y\mid x)$,
		the DA algorithm has Markov kernel
		\begin{equation} \nonumber
			K(x,\df x')=\int_{\Y}f_{Y\mid X}(\df y\mid x)f_{X\mid Y}(\df x'\mid y).
		\end{equation}
		It is reversible with respect to $f_X$. 
		Given the currents state~$x$, the latent variable $Y$ is drawn conditionally on $X = x$, and then~$X$ is updated conditionally on $Y$. 
		This is in fact the marginal chain of a two-block Gibbs sampler. We focus on two special cases.
		\begin{itemize}
			\item \namedlabel{ex:proximal}{\textsf{Alg-Proximal}}:
			If $\X=\Y=\mathbb{R}^d$ and $f_{Y\mid X}(\cdot\mid x)=\mathrm{N}_d(x,\eta I_d)$ for $\eta>0$, the DA algorithm defines a type of proximal sampler.
			\item \namedlabel{ex:gibbs}{\textsf{Alg-Gibbs}}:
			Suppose that $\X=\prod_{i=1}^B\X_i$ is a product of Polish spaces. The random-scan Gibbs sampler selects $J$ uniformly from $[B]:=\{1,\ldots,B\}$ and replaces $X_{J}$ by a draw from its conditional distribution given $X_{-J}$, retaining all other blocks. Here $X_{-j}$ denotes all blocks except $X_{j}$. This is a DA algorithm with
			\[
			\Y=\bigsqcup_{j=1}^B\{j\}\times\X_{-j},
			\qquad Y=(J,X_{-J}),
			\]
			where $\X_{-j}=\prod_{i\ne j}\X_i$, and $J$ is independent of $X$.
		\end{itemize}
		When exact sampling from $f_{X\mid Y}(\cdot\mid y)$ is impractical, it can be replaced by a transition $Q_y$ reversible with respect to that conditional distribution. The resulting hybrid DA algorithm remains reversible with respect to $f_X$. Comparisons between ideal and hybrid algorithms are available in \cite{andrieu2018uniform,qin2024spectral,qin2025spectral}; we will use them to analyze a hybrid proximal sampler and its choice of step size $\eta$.
	\end{itemize}
	
	We use mixture distributions as a surrogate for multi-modal or heavy-tailed target distributions.
	Assume that the target has the mixture representation
	\begin{equation} \label{eq:mixture}
		f_X(\cdot)=\int_{\Z}f_{X\mid Z}(\cdot\mid z)f_Z(\df z),
	\end{equation}
	where $\Z$ is Polish and $f_Z$ is a probability distribution.  
	To define $f_{X \mid Z}(\cdot \mid z)$, we use $\X=\mathbb{R}^d$ and a differentiable function $\xi$ such that $\exp(-\xi)$ is a probability density and
	\[
	\begin{aligned}
	\frac a2\|x'-x\|^2
	&\leq \xi(x')-\xi(x)-\nabla\xi(x)^\top(x'-x) \leq\frac L2\|x'-x\|^2,\qquad x,x'\in\mathbb R^d,
	\end{aligned}
	\]
	for constants $0<a\leq L<\infty$.
	Thus, $\xi$ is $a$-strongly convex and has an $L$-Lipschitz gradient.
	
	\begin{itemize}
		\item \namedlabel{ex:location}{\textsf{Location-mixture}}:
		For a measurable map $m:\Z\to\mathbb{R}^d$, the component density is
		$Df_{X\mid Z}(x\mid z)=\exp[-\xi(x-m(z))]$.
		Separated component locations allow multimodal targets.
		\item \namedlabel{ex:scale}{\textsf{Scale-mixture}}:
		Here $\Z=(0,\infty)$ and
		$Df_{X\mid Z}(x\mid z)=z^{-d/2}\exp[-\xi(x/\sqrt z)]$.
		We additionally assume that $\xi$ is twice continuously differentiable and that $\exp(-\xi)$ has mean zero.
		A heavy-tailed mixing distribution can produce a heavy-tailed target even though every component is strongly log-concave.
	\end{itemize}
	
	We analyze~$K$ by decomposing its dynamics into local and global levels.
	
	For each $z$, let $K_z$ be a local, or component Markov kernel reversible with respect to the component distribution $f_{X\mid Z}(\cdot\mid z)$. For \ref{ex:MH}, we let $K_z$ correspond to a Metropolis-Hastings algorithm targeting $f_{X \mid Z}$ with the same proposal as~$K$:
	\begin{equation} \label{eq:Kz-MH}
		K_z(x,\mathsf{A})=\int_{\mathsf{A}} \operatorname{acc}_z(x,x') \, q(x,x')\,\df x'
		+\left[1-\int_{\X} \operatorname{acc}_z(x,x') \, q(x,x')\,\df x'\right]\ind_{\mathsf{A}}(x),
	\end{equation}
	where
	\[
	\operatorname{acc}_z(x,x')=1\wedge\frac{Df_{X \mid Z}(x' \mid z)q(x',x)}{Df_{X \mid Z}(x \mid z)q(x,x')}.
	\]
	For the DA algorithm, we define $K_z$ to be a localized DA algorithm targeting the component distribution, based on the joint distribution $f_{X \mid Z}(\df x \mid z) \, f_{Y \mid X}(\df y \mid x)$:
	\begin{equation} \label{eq:Kz-DA}
		K_z(x, \df x') = \int_{\Y} f_{Y \mid X,Z}(\df y \mid x,z) \, f_{X \mid Y,Z}(\df x' \mid y,z) = \int_{\Y} f_{Y \mid X}(\df y \mid x) \, f_{X \mid Y,Z}(\df x' \mid y,z).
	\end{equation}
	For $z \in \Z$, $K_z$ describes the movement within a mixture component driven by~$K$.
	
	Information about the global movement across components driven by~$K$ is described by the following two Markov kernels on the mixture index space~$\Z$:
	\[
	\begin{aligned}
		\bar K_0(z,\df z')
		&=\int_{\X}f_{X\mid Z}(\df x\mid z) \, f_{Z\mid X}(\df z'\mid x),\\
		\bar K_1(z,\df z')
		&=\int_{\X^2}f_{X\mid Z}(\df x\mid z) \, K(x,\df x') \, f_{Z\mid X}(\df z'\mid x').
	\end{aligned}
	\]
	The first kernel measures the overlap between the component distributions $f_{X \mid Z}(\cdot \mid z)$. The second also incorporates one transition of~$K$ and can therefore capture movement not visible from overlap alone. Both kernels are reversible with respect to $f_Z$. 
	Convergence-related bounds for these kernels can be easier to obtain than bounds for $K$, especially when~$\Z$ has low dimension or a simple adjacency structure.
	
	Our main technical result, Theorem~\ref{thm:main}, combines weak Poincar\'{e} inequalities for $K_z$ with an inequality for either $\bar{K}_0$ or $\bar{K}_1$.
	It gives mixing-time bounds that can be applied to \ref{ex:MH}, \ref{ex:proximal}, and \ref{ex:gibbs} when $f_X$ is \ref{ex:location} or \ref{ex:scale}.
	Below, we give an illustration of the bounds in two concrete settings.
	
	\subsection{A glimpse of mixing-time results}

	Pick a total-variation tolerance $\epsilon>0$ and an initial law $\mu\ll f_X$ with
	\[
	\|\mu-f_X\|_{f_X}
	=\left[\int_{\X}\left(\frac{\df\mu}{\df f_X}-1\right)^2\df f_X\right]^{1/2}<\infty.
	\]
	Write $t(\epsilon;\mu,K)$ for the mixing time of the algorithm under consideration, defined to be the first time the chain reaches within $\epsilon$ of $f_X$ in total variation. 
	For illustration, we present below mixing time bounds in two special settings.
	In the first setting, we consider a multi-modal target, written into a location mixture, where a finite number of locations are connected along a path, but the path can be long.
	In the second setting, we consider a heavy-tailed target, written into a scale mixture, where we use an inverse gamma distribution for $f_Z$ to control the heaviness of the tails.
	
	\emph{Location mixtures with a path structure.}
	In \ref{ex:location}, assume that $\Z = [p] :=\{1,\dots,p\}$, $f_Z(z) = 1/p$ for $z \in [p]$, and $\|m(z+1)-m(z)\|\leq\kappa$ for $z\in[p-1]$. 
	For \ref{ex:MH} with $\sigma=\varsigma/\sqrt{Ld}$ and fixed $\varsigma>0$, \ref{ex:proximal} with step size $\eta$, and \ref{ex:gibbs} with $B$ blocks, Corollaries~\ref{cor:mixing-MH-location}, \ref{cor:mixing-proximal-location}, and~\ref{cor:mixing-gibbs-location} give
	\[
	t(\epsilon;\mu,K)=
	\begin{cases}
		O\!\left((L/a)d\,p^3e^{L\kappa^2/4}\log(\|\mu-f_X\|_{f_X}/\epsilon)\right),&\ref{ex:MH},\\
		O\!\left([1+(a\eta)^{-1}]p^3e^{L\kappa^2/[4(1+L\eta)]}\log(\|\mu-f_X\|_{f_X}/\epsilon)\right),&\ref{ex:proximal},\\
		O\!\left((L/a)B\,p^3e^{L\kappa^2/4}\log(\|\mu-f_X\|_{f_X}/\epsilon)\right),&\ref{ex:gibbs}.
	\end{cases}
	\]
	as $L/a$, $d$, $p$, $\|\mu-f_X\|_{f_X}/\epsilon$, $(a\eta)^{-1}$ (in \ref{ex:proximal}), $L \kappa^2$ (in \ref{ex:MH} and \ref{ex:gibbs}), $L\kappa^2/(1+L \eta)$ (in \ref{ex:proximal}), $B$ (in \ref{ex:gibbs}) tend to infinity.
	The exponential factor depends on neighboring separation $L \kappa^2$, while the number of components~$p$ enters polynomially. 
	These results remain informative as the mixture spreads over an increasingly large region, provided neighboring separation is uniformly controlled.
	
	Without any adjacency structure, if we assume $\sup_{z,z'} \|m(z) - m(z')\| \leq R$ for some $R \in [0,\infty)$, then, as $R \to \infty$, the mixing-time bounds we provide in Section \ref{sec:applications} scale with $e^{LR^2/4}$, or $e^{LR^2/[4(1+L\eta)]}$ for \ref{ex:proximal}.
	We note that tempering algorithms can potentially achieve a much better polylogarithmic dependence on~$R$ \citep{zhou2025polynomial}.
	
	In general, we provide spectral gap bounds for \ref{ex:location} when $f_Z$ is a generic probability measure on a possibly continuous space~$\Z$.
	The neighboring separation is characterized by the spectral gap of the transition kernel of the form
	\[
	\breve{K}(z, \mathsf{A}) = \int_{\mathsf{A}} \psi(z,z') \, f_Z(\df z') + r(z) \, \ind_{\mathsf{A}}(z),
	\]
	where $\psi(z,z') = \exp[-L\|m(z) - m(z')\|^2/4]$ for $z \neq z'$ (with minor modifications in the analysis of \ref{ex:proximal}).
	This transition kernel is reversible with respect to $f_Z$, and serves as a surrogate for $\bar{K}_0$ or $\bar{K}_1$.
	Rapid mixing can be ensured if $\breve{K}$ has a sufficient spectral gap, which in turn is guaranteed if the distribution of $m(Z)$ is sufficiently well-connected at scale $L^{-1/2}$.
	Weak Poincar\'{e} inequalities can relax this condition by allowing a small amount of probability outside a well-connected region.

	\emph{Scale mixtures with inverse-gamma mixing density.}
	In \ref{ex:scale}, let $Z\sim\mathrm{InvGamma}(\theta/2,\theta/2)$, where $\theta>0$.
	A smaller $\theta$ gives heavier tails. 
	Indeed, if $f_{X \mid Z}(\cdot \mid z)$ is the $\mathrm{N}_d(0, zI_d)$ distribution, then $f_X$ is the multivariate Student's $t$ distribution with degrees of freedom~$\theta$.
	In fact, one can show that, as $r \to \infty$, $\mathrm{Pr}(\|X\| \geq r)$ is of order $r^{-\theta}$.
	For \ref{ex:MH} with $\sigma=\varsigma/\sqrt{Ld}$ and fixed $0<\varsigma^2<1/2$, and for \ref{ex:proximal}, Corollaries~\ref{cor:mixing-MH-scale},~\ref{cor:mixing-proximal-scale}, and~\ref{cor:mixing-gibbs-scale} give
	\[
	\begin{aligned}
		t(\epsilon;\mu,K)=
		\begin{cases}
			O\!\left(
			(L/a)(1+L/a)d^2 \theta^{-2}
			\left[ 1024(1+L/a)d \, \theta^{-2} \|\mu-f_X\|_{f_X}^2/\epsilon^2 \right]^{2/\theta}
			\right), & \ref{ex:MH},\\
			O\!\left((1+L/a)d (a\eta)^{-1} \theta^{-1}
			\left[1024(1+L/a)d \,
			\theta^{-2} \|\mu-f_X\|_{f_X}^2/\epsilon^2 \right]^{2/\theta}
			\right),&\ref{ex:proximal},\\
			O\!\left( (L/a)^2 B d \theta^{-2} \log(\|\mu-f_X\|_{f_X}/\epsilon)\right),&\ref{ex:gibbs}.
		\end{cases}
	\end{aligned}
	\]
	as $L/a,d,\|\mu-f_X\|_{f_X}/\epsilon$,  $\theta^{-1}$, $(a\eta)^{-1}$ (in \ref{ex:proximal}), $B$ (in \ref{ex:gibbs}) tend to infinity.
	In the first two algorithms, $\|\mu - f_X\|_{f_X}/\epsilon$ enters polynomially rather than logarithmically, so the bounds guarantee rapid convergence only when this ratio scales polynomially with $L/a$ and~$d$.
	This is because the parameters $\sigma$ and $\eta$ do not adapt to the value of $z$, so $K_z$ can mix poorly for extreme values of~$z$, leading to deteriorated mixing-time bounds for~$K$.
	On the other hand, for the Gibbs sampler, the spectral gap of $K_z$ is invariant under isotropic scaling, resulting in better dependence on warmness and error tolerance.
	A simulation study in Section \ref{ssec:scale-simulation} shows that \ref{ex:gibbs} appears to converge faster than \ref{ex:MH} and \ref{ex:proximal} when $\theta$ is small.
	
	More generally, we consider a generic distribution $f_Z$ on $(0,\infty)$.
	Since the overlap of scale components depend on relative scales, it makes sense to consider the distribution of $\log Z$, denoted by $f_{\log Z}$.
	Let $Df_{\log Z}(u)$ be the density of its continuous part.
	The spectral gap of $\bar{K}_0$ is shown to be of at least a positive constant multiple of $1 \wedge [\mathrm{Is}(f_{\log Z})^2 (1+L/a)^{-1} d^{-1}]$, where
	\begin{equation} \label{eq:is}
	\mathrm{Is}(f_{\log Z}) = \mathrm{ess \, inf}_{\ell<u<r} \, \frac{Df_{\log Z}(u)}{\min\{f_{\log Z}((-\infty,u]),\, f_{\log Z}((u, \infty))\}}
	\end{equation}
	is an isoperimetric constant quantifying the connectivity of $f_{\log Z}$.
	Here $\ell=\inf\{u:f_{\log Z}((-\infty,u])>0\}$, $r=\sup\{u:f_{\log Z}((-\infty,u])<1\}$, and the essential infimum is with respect to Lebesgue measure.
	For the aforementioned inverse gamma distribution, $\mathrm{Is}(f_{\log Z}) \geq \theta/(4\sqrt{\theta+2})$, giving the $\theta^{-2}$ dependence in the bound for the Gibbs sampler.
	For the other two samplers,  additional price is paid because $K_z$ does not mix well at every scale, and the probability mass of the bad scales must be taken into account.

	\subsection{Related work}

	Decomposition methods relate convergence within subsets or components to movement between them \citep{caracciolo1992two,madras2002markov}. The soft decomposition of \cite{ge2018simulated} provides a comprehensive framework involving spectral gap decomposition; it is developed principally for simulated tempering. 
	Both the spectral-gap and weak Poincar\'{e} decomposition bounds developed here fit within the general operator framework of \citet{qin2025spectral}, whose Section~6.2 combines weak inequalities for component dynamics and an idealized global chain. A related weak decomposition appears in \citet[Lemma~4.11]{huang2025weak}, under a positive global spectral-gap assumption.
	See also \cite{atchade2021approximate}, who considered the decomposition of approximate spectral gaps.
	
	\citet{gao2026weak} develop weak Poincar\'{e} comparisons for deterministic-scan Metropolis-within-Gibbs samplers with two blocks, obtaining explicit subgeometric convergence bounds for nonreversible chains. We adapt their comparison between conjugate operators to transfer weak Poincar\'{e} bounds between the mixture-index kernels and their counterparts on the target space.
	
	Our contribution is to specialize these principles to mixture targets through explicit comparisons involving the overlap kernel $\bar K_0$ and the algorithm-dependent index kernel $\bar K_1$, and to turn these comparisons into quantitative mixing-time bounds for the samplers studied here. The resulting location-mixture bounds capture neighboring-component separation along paths, while the scale-mixture bounds quantify the dependence on tail heaviness and initialization.

	Proximal samplers for mixture targets is studied in \citet[Proposition~2.4]{CLY2026}. They establish the spectral gap decomposition
	\[
	\mathrm{Gap}(K)\geq\delta\mathrm{Gap}(\bar K_1)
	\qquad\text{when }\inf_z\mathrm{Gap}(K_z)\geq\delta>0.
	\]
	We provide this bound in the more general DA setting in Corollary~\ref{cor:uniform-gap}, and we also give more general versions that involve weak Poincar\'{e} inequalities.
	
	Weak Poincar\'{e} inequalities and their conversion into subgeometric mixing bounds are developed systematically in \cite{andrieu2026weak}.
	The Student’s $t$ targets studied by \citet{andrieu2026weak} are the Gaussian special case of our inverse-gamma scale mixtures. For fixed $\theta>2$, their mixing-time bound for \ref{ex:MH}, gives an $O_\theta(d^3)$ dependence. 
	Our argument yields $O_\theta(d^{2+2/\theta})$, using a different proposal scaling. More broadly, our analysis covers non-Gaussian strongly log-concave components and explicitly tracks the deterioration of the bounds as $\theta\downarrow0$.
	
	For \ref{ex:gibbs}, \citet[Corollary~3.7]{ascolani2024entropy} provides the estimates for $\mathrm{Gap}(K_z)$ necessary for our analysis.
	
	The analysis of hybrid proximal samplers is based on comparison results from \cite{andrieu2018uniform,qin2024spectral,qin2025spectral}.
	
%	For DA algorithms, Appendix~\ref{app:approx} also gives a perturbation comparison, based on \cite{ascolani2024dimension}, between chains with a common augmentation conditional and nearby target distributions. This provides a way to transfer mixture-based convergence bounds to suitable approximations of other targets.

	\citet{chen2021dimension} establish Poincaré and log-Sobolev inequalities for mixtures from uniform component inequalities and moment bounds on pairwise $\chi^2$ divergences. In particular, they obtain dimension-free log-Sobolev bounds for Gaussian location mixtures whose mixing distributions have bounded support.
	Although not pursued here, it may be possible to utilize these functional inequalities to study Markov chains targeting mixture distributions.

	\subsection{Article Organization}
	
	The rest of the article is organized as follows.
	In Section \ref{sec:Poincare}, we summarize some facts about Poincar\'{e} inequalities for Markov chains and their implications.
	Section \ref{sec:decomposition} gives the general decomposition result that separates within- and between-component dynamics.
	Section \ref{sec:Kbar} gives some general results regarding the spectral gaps of the global transition kernels $\bar{K}_0$ and $\bar{K}_1$.
	The decomposition results are applied to the three algorithms targeting location and scale mixtures in Section \ref{sec:applications}.

	\section{Preliminaries: Poincar\'{e} Inequalities} \label{sec:Poincare}
	
	Let $(\Omega, \mathcal{F}, \varpi)$ be a probability space.
	Let $T: \Omega \times \mathcal{F} \to [0,1]$ be an Mtk that is reversible with respect to~$\varpi$.
	For $g: \Omega \to \mathbb{R}$ such that $\|g\|_{\varpi}^2 := \int_{\Omega} g(\omega)^2 \, \varpi(\df \omega) < \infty$, define the Dirichlet form
	\[
	\mathcal{E}_T(g) = \frac{1}{2} \int_{\Omega^2} \varpi(\df \omega) \, T(\omega, \df \omega') [g(\omega) - g(\omega')]^2 = \|g\|_{\varpi}^2 - \int_{\Omega^2} \varpi(\df \omega) \, T(\omega, \df \omega') \, g(\omega) g(\omega').
	\]
	Let $L_0^2(\varpi)$ be the space of measurable functions $g: \Omega \to \mathbb{R}$ such that $\varpi g := \int g \, \df \varpi = 0$ and $\|g\|_{\varpi}^2 < \infty$.
	Given $s \geq 0$ and $b \geq 0$, we say $T$ satisfies an $(s,b)$-Poincar\'{e} inequality if, for every measurable function $g: \Omega \to \mathbb{R}$ such that $g \in L_0^2(\varpi)$,
	\[
	\|g\|_{\varpi}^2 \leq s \, \mathcal{E}_T(g) + b \|g\|_{\mathrm{osc}}^2,
	\]
	where $\|g\|_{\mathrm{osc}}$ is the difference between the essential supremum and essential infimum of~$g$ defined in terms of $\varpi$.
	In the special case that $T$ satisfies an $(s,0)$-Poincar\'{e} inequality, i.e., $\|g\|_{\varpi}^2 \leq s \, \mathcal{E}_T(g)$ for every $g\in L_0^2(\varpi)$,
	we say $T$ admits a spectral gap of size $1/s$, and write $\mathrm{Gap}(T) \geq 1/s$.
	For the degenerate case $L_0^2(\varpi)=\{0\}$, we use the convention $\mathrm{Gap}(T)=1$ and interpret spectral-gap lower bounds as being at most one.
	
	\begin{remark}
		By Popoviciu's inequality, $\|g\|_{\varpi}^2 \leq \|g\|_{\mathrm{osc}}^2/4$ for $g \in L_0^2(\varpi)$.
		Thus, for any $s \geq 0$, $T$ always satisfies an $(s,b)$-Poincar\'{e} inequality if $b \geq 1/4$.
	\end{remark}

	The above Poincar\'{e} inequality has some important consequences regarding the convergence properties of~$T$.
	Before giving the bounds, we introduce some additional concepts and notations.
	
	In what follows, for a distribution $\mu: \mathcal{F} \to [0,1]$, write $\mu T^t(\cdot) = \int_{\Omega} \mu(\df \omega) \, T^t(\omega, \cdot)$, where $T^0(\omega, \mathsf{A}) = \ind_{\mathsf{A}}(\omega)$ for $\omega \in \Omega$ and $\mathsf{A} \in \mathcal{F}$, and $T^t(\omega, \cdot) = \int_{\Omega} T(\omega, \df \omega') \, T^{t-1}(\omega', \cdot)$ for $\omega \in \Omega$ and $t \geq 1$.
	
	Denote by $L_*^2(\varpi)$ the set of probability measures~$\mu$ such that $\mu \ll \varpi$ and that $\|\df \mu / \df \varpi\|_{\varpi} < \infty$.
	For two distributions $\mu$ and $\nu$ in $L_*^2(\varpi)$, let $\|\mu - \nu\|_{\varpi}$ be their $L^2(\varpi)$ distance. When $\nu=\varpi$, its square is the $\chi^2$ divergence.
	That is, $\|\mu - \nu\|_{\varpi} = \|\df \mu/\df\varpi - \df \nu / \df \varpi\|_{\varpi} = \sup |\mu g - \nu g|$, where the supremum is over functions $g \in L_0^2(\varpi)$ such that $\|g\|_{\varpi} = 1$.
	Also, let $\mathrm{TV}(\mu, \nu)$ be the total variation distance between~$\mu$ and~$\nu$.
	It is well-known that $\mathrm{TV}(\mu, \nu) \leq (1/2) \|\mu - \nu\|_{\varpi}$.

	The following result is established in Appendix~\ref{app:convergence-Poincare}.
	
	\begin{proposition} \label{pro:convergence-Poincare}
		Suppose that $T$ is positive semi-definite in the sense that $$\int_{\Omega^2} \varpi(\df \omega) \, T(\omega, \df \omega') \, g(\omega) \, g(\omega') \geq 0$$ for $g \in L_0^2(\varpi)$, and that it satisfies an $(s,b)$-Poincar\'{e} inequality for some $b \geq 0$ and $s \geq 1$.
		Then, for any distribution $\mu \in L_*^2(\varpi)$ and positive integer $t$,
		\[
		\begin{aligned}
			\mathrm{TV}(\mu T^t, \varpi)^2 \leq \|\mu - \varpi\|_{\varpi}^2 \left\{ \frac{1}{4} \left(1 - \frac{1}{s} \right)^t + b\left[1-\left(1-\frac{1}{s}\right)^t\right] \right\}.
		\end{aligned}
		\]
	\end{proposition}
	
	\begin{remark}
		It is well-known that the Mtk of a DA algorithm is positive semi-definite.
		Slightly less well-known is the fact that the Mtk associated with \ref{ex:MH} is positive semi-definite as well \citep{baxendale2005renewal}.
		One can also obtain a positive semi-definite Mtk from an arbitrary~$T$ by considering a lazy version of it, given by $T_{\mathrm{lazy}}(\omega, \mathsf{A}) = [T(\omega, \mathsf{A}) + \ind_{\mathsf{A}}(\omega)]/2$ for $\omega \in \Omega$ and $\mathsf{A} \in \mathcal{F}$.
	\end{remark}
	
	For $\epsilon > 0$, define the mixing time $t(\epsilon; \mu, T)$ to be the first time that $\mathrm{TV}(\mu T^t, \varpi)$ reaches~$\epsilon$ or below.
	Then Corollary~\ref{cor:mixing-sb} below follows from Proposition~\ref{pro:convergence-Poincare}.
	\begin{corollary} \label{cor:mixing-sb}
		Assume that $T$ is positive semi-definite, and it satisfies an $(s,b)$-Poincar\'{e} inequality, where $s > 1$, and $\sqrt{b} < \epsilon / \|\mu-\varpi\|_{\varpi} < 1/2$.
		Then, for any distribution $\mu \in L_*^2(\varpi)$,
		\[
		t(\epsilon; \mu, T) \leq \left\lceil\frac{\log \left[ (\epsilon^2/\|\mu - \varpi\|_{\varpi}^2 - b)/(1/4-b) \right]}{\log (1-1/s)}\right\rceil.
		\]
	\end{corollary}

	When a positive-definite or lazy $T$ satisfies an $(s,0)$-Poincar\'{e} inequality for some $s \geq 1$, i.e., $\mathrm{Gap}(T) \geq 1/s$, $\mu T^t$ approaches $\varpi$ at the geometric rate $(1-1/s)^{t/2}$.
	In fact, a spectral gap of size $1/s$ can guarantee a faster convergence rate of $(1-1/s)^t$ \citep{roberts1997geometric}, although this is not reflected in Proposition~\ref{pro:convergence-Poincare}.
	When $b \neq 0$, Proposition \ref{pro:convergence-Poincare} does not provide an upper bound on $\mathrm{TV}(\mu T^t,\varpi)$ that goes to~0 as $t \to \infty$.
	This can be overcome if one can establish a collection of $(s,b)$-Poincar\'{e} inequalities.
	
	Let $\beta: (0,\infty) \to [0,\infty)$ be a measurable function.
	A weak Poincar\'{e} inequality holds for~$T$ with rate function~$\beta$ if, for every $s > 0$, $T$ satisfies an $(s, \beta(s))$-Poincar\'{e} inequality.

	\begin{proposition} \label{pro:Poincare-weak} \citep{andrieu2026weak}
		Suppose that $T$ is positive semi-definite, and that it satisfies a weak Poincar\'{e} inequality with rate function $\beta: (0,\infty) \to [0,\infty)$.
		Assume that $\beta$ is non-increasing, and $\beta(s) \to 0$ as $s \to \infty$.
		Then, for any distribution $\mu \in L_*^2(\varpi)$ and positive integer~$t$,
		\[
		\mathrm{TV}(\mu T^t, \varpi)^2 \leq \|\mu-\varpi\|_{\varpi}^2 \, G_{\beta}^{-1}(t), 
		\]
		where
		\[
		G_{\beta}(v) = \int_v^{1/4} \frac{1}{H_{\beta}(u)} \, \df u , \quad H_{\beta}(u) = \sup_{s > 0} \frac{u - \beta(s)}{s}.
		\]
		Moreover, as $t \to \infty$, the upper bound goes to~0.
	\end{proposition}
	\begin{proof}
		By Lemma~\ref{lem:spectral} in Appendix~\ref{app:operator}, for any $g \in L_0^2(\varpi)$, $\langle g, T g \rangle_{\varpi} \geq \langle g, T^2 g \rangle_{\varpi}$.
		Thus, $\mathcal{E}_T(g) \leq \mathcal{E}_{T^2}(g)$.
		In particular,~$T^2$ satisfies a weak Poincar\'{e} inequality with rate $\beta$.
		The desired result then follows from Theorem~8 and Proposition~24 of \cite{andrieu2026weak}.
	\end{proof}
	
	Weak Poincar\'{e} inequalities are often used to establish sub-geometric convergence for Markov chains.
	Roughly speaking, the faster $\beta(s)$ converges to~0 as $s \to \infty$, the slower $H_{\beta}(u)$ converges to~0 as $u \to 0$, and the faster $G_{\beta}^{-1}(t)$ converges to~0 as $t \to \infty$.
	
	\begin{proposition} \label{pro:polynomial}
		In Proposition~\ref{pro:Poincare-weak}, if $\beta(s) = c s^{-\ell}$ for some $\ell > 0$ and $c > 0$, then $H_{\beta}(u) = \ell (\ell+1)^{-( \ell +1)/ \ell } c^{-1/ \ell } u^{( \ell +1)/ \ell }$, and $G_{\beta}^{-1}(t) \leq c ( \ell +1)^{ \ell +1} t^{- \ell }$.
	\end{proposition}
	\begin{proof}
		The supremum defining $H_{\beta}(u)$ is attained at
		$s=[( \ell +1)c/u]^{1/ \ell }$, which gives the stated formula. Consequently,
		\[
		G_\beta(v)=( \ell +1)^{1+1/ \ell }c^{1/ \ell }
		\bigl(v^{-1/ \ell }-4^{1/ \ell }\bigr)
		\leq ( \ell +1)^{1+1/ \ell }c^{1/ \ell }v^{-1/ \ell }.
		\]
		Inverting this inequality gives the bound on $G_\beta^{-1}(t)$.
	\end{proof}
	
	\begin{corollary} \label{cor:mixing-polynomial}
		Suppose that $T$ is positive semi-definite and satisfies a weak
		Poincar\'{e} inequality with $\beta(s)\leq c s^{-\ell}$ for some
		$c>0$ and $\ell>0$. For any $\mu\in L_*^2(\varpi)$ and $\epsilon>0$,
		\[
		t(\epsilon;\mu,T)\leq
		\left\lceil\left(
		\frac{c(\ell+1)^{\ell+1}\|\mu-\varpi\|_{\varpi}^2}{\epsilon^2}
		\right)^{1/\ell}\right\rceil.
		\]
		%In particular, for fixed $\ell$ and $\epsilon$, as
		%$c^{1/\ell}\|\mu-\varpi\|_{\varpi}^{2/\ell}\to\infty$,
		%\[
		%t(\epsilon;\mu,T)
		%=O\bigl(c^{1/\ell}\|\mu-\varpi\|_{\varpi}^{2/\ell}\bigr),
		%\]
		%with an implied constant depending only on $\ell$ and $\epsilon$.
	\end{corollary}
	\begin{proof}
		Propositions~\ref{pro:Poincare-weak} and~\ref{pro:polynomial} give
		\[
		\mathrm{TV}(\mu T^t,\varpi)^2
		\leq c(\ell+1)^{\ell+1}\|\mu-\varpi\|_{\varpi}^2t^{-\ell},
		\qquad t\geq1.
		\]
		Solving for $t$ and rounding up proves the claim when $\mu\ne\varpi$;
		when $\mu=\varpi$, the mixing time is zero.
	\end{proof}
	
	For polynomial rates, Corollary~\ref{cor:mixing-polynomial} can yield a better mixing-time dependence compared to
	Corollary~\ref{cor:mixing-sb} with a single choice of $s$.

	\section{Markov Decomposition} \label{sec:decomposition}
	
	Recall that~$K$ is an Mtk that is reversible with respect to the target distribution $f_X$, and the target distribution has the mixture representation~\eqref{eq:mixture}.
	Suppose that, for $z \in \Z$, $K_z$ is an Mtk on~$\X$ that is reversible with respect to $f_{X \mid Z}(\cdot \mid z)$.
	We assume that $(z,x)\mapsto K_z(x,\mathsf A)$ is measurable for every $\mathsf A\in\mathcal B$. We work under the following technical condition:
	\begin{itemize}
		\item \namedlabel{cond:A}{(A)} 
		For $g: \X \to \mathbb{R}$ such that $\|g\|_{f_X}^2 < \infty$,
		\[
		\mathcal{E}_K(g) \geq \int_{\Z} f_Z(\df z) \, \mathcal{E}_{K_z}(g).
		\]
	\end{itemize}
	
	This condition has appeared in \cite{ge2018simulated} and subsequent works.
	By defining $K_z$ appropriately, it is possible to establish Condition \ref{cond:A} for a wide range of samplers \citep{lee2024convergence}.
	For instance, the following holds.
	
	\begin{proposition} \label{pro:condition-MH} \citep{lee2024convergence}
		Suppose that $K$ is (the Mtk associated with) the MH algorithm targeting $f_X$ with proposal density $q(\cdot,\cdot)$.
		For $z \in \Z$, let $K_z$ be the MH algorithm targeting $f_{X \mid Z}(\cdot \mid z)$ with the same proposal density, i.e., assume \eqref{eq:Kz-MH}.
		Then Condition~\ref{cond:A} holds.
	\end{proposition}
	
	To analyze the DA algorithm, we establish the following proposition in Appendix~\ref{app:condition-DA}.
	\begin{proposition} \label{pro:condition-DA}
		Suppose that $K$ is the DA algorithm targeting $f_X$ based on the joint distribution $f_X(\df x) \, f_{Y \mid X}(\df y \mid x)$.
		Let $K_z$ be the DA algorithm targeting $f_{X \mid Z}(\cdot \mid z)$ based on the joint distribution $f_{X \mid Z}(\df x \mid z) \, f_{Y \mid X}(\df y \mid x)$.
		That is, assume \eqref{eq:Kz-DA}.
		Then Condition \ref{cond:A} holds.
	\end{proposition}
	
	\begin{remark} \label{rem:condition-DA}
		In Proposition \ref{pro:condition-DA}, to simulate $K_z$, given the current state $x \in \X$, the latent variable~$y$ is drawn from $f_{Y \mid X}(\cdot \mid x)$ --- this is the same as in the original DA algorithm~$K$.
		If $K$ is \ref{ex:proximal} targeting $f_X$ associated with a Gaussian family $\mathrm{N}(x, \eta I_d)$, then $K_z$ would be a proximal sampler targeting $f_{X \mid Z}(\cdot \mid z)$ associated with the same Gaussian family.
		If $K$ is \ref{ex:gibbs} targeting a $f_X(\df x_1, \dots, \df x_B)$, then $K_z$ would be a random-scan Gibbs algorithm targeting $f_{X \mid Z}(\df x_1, \dots, \df x_B \mid z)$ under the same blocking scheme.
	\end{remark}

	Recall also the global Mtks on~$\Z$:
	\[
	\bar{K}_0(z,\df z') = \int_{\X} f_{X \mid Z}(\df x \mid z) \, f_{Z \mid X}(\df z' \mid x), \quad \bar{K}_1(z, \df z') = \int_{\X^2} f_{X \mid Z}(\df x \mid z) \, K(x, \df x') \, f_{Z \mid X}(\df z' \mid x') .
	\]
	
	In \cite{ge2018simulated}, the relationship between the spectral gaps of~$K$ and $K_z$ is examined, with a focus on simulated tempering algorithms.
	Theorem~\ref{thm:main} below can be seen as an extension of their work, combining their idea (which in turn stems from \cite{caracciolo1992two} and \cite{madras2002markov}) with the notion of $(s,b)$-Poincar\'{e} inequalities.
	An analogous local/global argument appears in Lemma 4.11 of \cite{huang2025weak}; that result allows the local chains to satisfy weak Poincar\'{e} inequalities, while assuming the global chain $\bar{K}_0$ admits a spectral gap. 
	Here, we allow both the local and global chains to satisfy only weak Poincar\'{e} inequalities, and we allow the use of $\bar{K}_1$ in our analysis.
	We provide the proof in Appendix~\ref{app:main}.
	
	\begin{theorem} \label{thm:main}
		Assume that Condition \ref{cond:A} holds.
		Suppose that there exists $s \geq 0$ such that, for $f_Z$-almost every $z \in \Z$, the Mtk $K_z$ satisfies an $(s,b_z)$-Poincar\'{e} inequality, where $b_z \geq 0$, and $z \mapsto b_z$ is measurable. We may take $b_z\leq1/4$ without loss of generality, by Popoviciu's inequality.
		Then each of the following statements holds:
		\begin{enumerate}
			\item [(i)] If $\bar{K}_0$ satisfies an $(s',b')$-Poincar\'{e} inequality for some non-negative numbers $s'$ and $b'$, then $K$ satisfies an $(s'',b'')$-Poincar\'{e} inequality, where $s'' = (s'+1) s$, and $b'' = b' + (s'+1) \int_{\Z} b_z \, f_Z(\df z) $.
			\item [(ii)] If $K$ is positive semi-definite and $\bar{K}_1$ satisfies an $(s',b')$-Poincar\'{e} inequality for some non-negative numbers $s'$ and $b'$, then $K$ satisfies an $(s'',b'')$-Poincar\'{e} inequality, where $s'' = (s'+1)(s+2)$, and $b''$ is as in (i).
		\end{enumerate}
	\end{theorem}

	\begin{remark}
		It is shown in Proposition~\ref{pro:converse} that, if $\bar{K}_0$ satisfies an $(s,b)$-Poincar\'{e} inequality, then $\bar{K}_1$ satisfies an $(s,b)$-Poincar\'{e} inequality with the same value of $(s,b)$.
		It stands to reason that Part (ii) of Theorem~\ref{thm:main} has the potential to yield a sharper Poincar\'{e} inequality for~$K$, compared to Part (i).
		On the other hand, $\bar{K}_0$ takes a simpler form compared to $\bar{K}_1$, so Part (i) is typically easier to apply.
		In most of the examples we consider, we apply Part (i) instead of (ii).
	\end{remark}

	In situations where $\bar{K}_0$ (or $\bar{K}_1$) and each $K_z$ admits a spectral gap, the following proposition, whose proof is given in Appendix~\ref{app:main-proposition}, is convenient.
	
	\begin{proposition} \label{pro:main}
		Assume that Condition \ref{cond:A} holds.
		Suppose that there is a measurable function $\varphi: \Z \to (0,\infty]$ such that $\mathrm{Gap}(K_z) \geq 1/\varphi(z)$ for $f_Z$-almost every $z \in \Z$.
		Fix $s \geq 0$, and let
		\[
		b_z = \begin{cases}
			1/4 & \text{if } \; \varphi(z) > s, \\
			0 & \text{otherwise}.
		\end{cases}
		\]
		Then each of the following statements holds:
		\begin{enumerate}
			\item [(i)] If $\mathrm{Gap}(\bar{K}_0) \geq 1/s'$ for some $s' \geq1$, then $K$ satisfies an $(s'',b'')$-Poincar\'{e} inequality, where $s'' = s's$, and $b'' = s'\int_{\Z} b_z \, f_Z(\df z)$.
			\item [(ii)] If $K$ is positive semi-definite and $\mathrm{Gap}(\bar{K}_1) \geq 1/s'$ for some $s' \geq1$, then $K$ satisfies an $(s'',b'')$-Poincar\'{e} inequality, where $s'' = s'(s+2)$, and $b''$ is as in (i).
		\end{enumerate}
	\end{proposition}

	We also present a corollary of Proposition \ref{pro:main} that gives the spectral gap decomposition property \citep{madras2002markov,qin2025spectral}.
	Appendix \ref{app:uniform-gap} contains the proof.

	\begin{corollary}\label{cor:uniform-gap}
		Assume Condition~\ref{cond:A}, and suppose that $\mathrm{Gap}(K_z)\geq\delta>0$ for $f_Z$-almost every $z$. Then
		\[
		\mathrm{Gap}(K)\geq\delta \mathrm{Gap}(\bar{K}_0).
		\]
		If, in addition, $K$ is positive semi-definite, then
		\[
		\mathrm{Gap}(K)\geq\frac{\delta\mathrm{Gap}(\bar K_1)}{1+2\delta}.
		\]
		If $K$ and $K_z$ are the DA kernels of Proposition~\ref{pro:condition-DA}, then we have the sharper bound
		\begin{equation} \nonumber
			\mathrm{Gap}(K)\geq\delta\mathrm{Gap}(\bar K_1)
		\end{equation}
	\end{corollary}

	Additionally, it is easy to use Theorem \ref{thm:main} to combine weak Poincar\'{e} inequalities for $K_z$ and $\bar{K}_t$, $t \in \{0,1\}$, and derive a weak Poincar\'{e} inequality for~$K$.
	
	\begin{corollary} \label{cor:weak}
		Assume that Condition \ref{cond:A} holds.
		Suppose that, for $f_Z$-almost every $z \in \Z$, $K_z$ satisfies a weak Poincar\'{e} inequality with rate function $\alpha_z: (0,\infty) \to [0,\infty)$.
		The function $(z,s) \mapsto \alpha_z(s)$ is assumed to be measurable.
		Then each of the following statements holds:
		\begin{enumerate}
			\item[(i)] If $\bar{K}_0$ satisfies a weak Poincar\'{e} inequality with rate function $\beta: (0,\infty) \to [0,\infty)$, then $K$ satisfies a weak Poincar\'{e} inequality with any measurable rate function bounded below by
			\[
			s'' \mapsto \gamma(s'') = \begin{cases}
				\inf_{s,s'} \left[ \beta(s') + (s'+1) \int_{\Z} \alpha_z(s) \, f_Z(\df z) \right] & \text{if } s'' > 1, \\
				1/4 & \text{otherwise},
			\end{cases}
			\]
			where the infimum is taken over the set of positive numbers~$s$ and~$s'$ such that $s'' = (s'+1)s$.
			\item[(ii)] If $K$ is positive semi-definite and $\bar{K}_1$ satisfies a weak Poincar\'{e} inequality with rate function $\beta: (0,\infty) \to [0,\infty)$, then $K$ satisfies a weak Poincar\'{e} inequality with any measurable rate function bounded below by
			\[
			s'' \mapsto \gamma(s'') = \begin{cases}
				\inf_{s,s'} \left[ \beta(s') + (s'+1) \int_{\Z} \alpha_z(s) \, f_Z(\df z) \right] & \text{if } s'' > 2, \\
				1/4 & \text{otherwise},
			\end{cases}
			\]
			where the infimum is taken over the set of positive numbers~$s$ and~$s'$ such that $s'' = (s'+1)(s+2)$.
		\end{enumerate}
	\end{corollary}
	
	To end this subsection, we state a result that is complementary to Theorem~\ref{thm:main}, the proof of which can be found in Appendix~\ref{app:converse}.
	
	\begin{proposition} \label{pro:converse}
		If $\bar{K}_0$ satisfies an $(s,b)$-Poincar\'{e} inequality for some $s \geq 0$ and $b \geq 0$, then $\bar{K}_1$ satisfies a Poincar\'{e} inequality with the same value of $(s,b)$.
		If $K$ satisfies an $(s,b)$-Poincar\'{e} inequality for some $s \geq 1$ and $b \geq 0$, then~$\bar{K}_1$ satisfies a Poincar\'{e} inequality with the same value of $(s,b)$.
	\end{proposition}

	\section{Spectral Gaps of \texorpdfstring{$\bar{K}_0$ and $\bar{K}_1$}{the Mixture-Index Kernels}} \label{sec:Kbar}
	
	\subsection{General results}
	
	To apply Theorem~\ref{thm:main}, we would need to establish Poincar\'{e} inequalities for the $K_z$'s and for $\bar{K}_0$ or $\bar{K}_1$.
	The local Mtk $K_z$ will be examined on a case-by-case basis, when~$K$ corresponds to one of \ref{ex:MH}, \ref{ex:proximal}, and \ref{ex:gibbs}, and when $f_X$ is \ref{ex:location} or \ref{ex:scale}.
	For now, we derive some general properties for $\bar{K}_0$ and $\bar{K}_1$.
	In particular, we focus on $\bar{K}_0$, which is fully determined by $f_Z$ and $f_{X \mid Z}$.
	
	Recall that $\bar{K}_0(z, \df z') = \int_{\X} f_{X \mid Z}(\df x \mid z) \, f_{Z \mid X}(\df z' \mid x)$.
	Thus, $\bar{K}_0$ can be seen as a DA algorithm associated with the random element pair $(Z,X)$, which has joint distribution $f_Z(\df z) \, f_{X \mid Z}(\df x \mid z)$.
	Thus, it has all the properties of a generic DA algorithm, listed in Appendix~\ref{app:da}.
	In particular, writing $\rho_{\max}(X,Z)$ for the maximal correlation between $X$ and $Z$,
	\begin{equation} \nonumber
		\mathrm{Gap}(\bar{K}_0) = 1-\rho_{\max}(X,Z)^2,\qquad
		\rho_{\max}(X,Z)=\sup_{g,h} \frac{\mathrm{Cov} [g(X), h(Z)] }{\sqrt{\mathrm{Var}[g(X)] \, \mathrm{Var}[h(Z)] }},
	\end{equation}
	where the supremum is over measurable functions with finite, strictly positive variances.
	See \cite{liu1994covariance}.
	Set the maximal correlation to zero if either centered $L^2$ space is trivial.
	Thus, the less dependence there is between~$X$ and~$Z$, the faster $\bar{K}_0$ converges.

	The conditional $f_{Z \mid X}(\cdot \mid x)$ typically has a somewhat complicated form, hindering the analysis of $\bar{K}_0$.
	To circumvent this, we use the following proposition, adopted from \cite{zhou2025polynomial}.
	Results of similar flavors can be found in, e.g., \cite{chen2021dimension,ge2018simulated}.
	
	\begin{proposition} \label{pro:peskun-DA-0}
		Suppose that, with respect to some $\sigma$-finite underlying measure on $(\X, \mathcal{B})$, for $f_Z$-almost every $z \in \Z$, $f_{X \mid Z}(\cdot \mid z)$ has density function $Df_{X \mid Z}(\cdot \mid z)$.
		Let $\breve{K}_0(z, \df z')$ be an Mtk that is reversible with respect to $f_Z$.
		Suppose that, for $z \in \Z$ and measurable $\mathsf{A} \subset \Z$,
		\[
		\breve{K}_0(z, \mathsf{A} \setminus \{z\}) \leq \int_{\mathsf{A}} \left[ \int_{\X} \sqrt{Df_{X \mid Z}(x \mid z) \, Df_{X \mid Z}(x \mid z')} \, \df x \right]^2 \, f_Z(\df z').
		\]
		Then, $\bar{K}_0$ satisfies an $(s,b)$-Poincar\'{e} inequality for some $s,b \geq 0$ whenever $\breve{K}_0$ satisfies an $(s,b)$-Poincar\'{e} inequality with the same value of $(s,b)$.
	\end{proposition}
	
	\begin{proof}
		By Lemma~\ref{lem:da-elements} in Appendix~\ref{app:da}, for $g \in L_0^2(f_Z)$, it holds that
		\[
		\begin{aligned}
			\mathcal{E}_{\bar{K}_0}(g) &= \frac{1}{2} \int_{\Z^2} f_Z(\df z) \, \bar{K}_0(z, \df z') \, [g(z) - g(z')]^2 \\
			&\geq \frac{1}{2} \int_{\Z^2} f_Z(\df z) \, f_Z(\df z') \, \left[ \int_{\X} \sqrt{Df_{X \mid Z}(x \mid z) \, Df_{X \mid Z}(x \mid z')} \, \df x \right]^2  [g(z) - g(z')]^2 \\
			&\geq \frac{1}{2} \int_{\Z^2} f_Z(\df z) \, \ind_{\Z \setminus \{z\}}(z') \, \breve{K}_0(z, \df z') \, [g(z) - g(z')]^2 = \mathcal{E}_{\breve{K}_0}(g).
		\end{aligned}
		\]
		This gives the desired result.
	\end{proof}
	
	When $K$ is the DA algorithm based on the random variable pair $(X,Y)$, one can establish a similar result for $\bar{K}_1$.
	Indeed, in this case, we can let $(Z,X,Y)$ have joint distribution $f_Z(\df z) \, f_{X \mid Z}(\df x \mid z) \, f_{Y \mid X}(\df y \mid x)$.
	For $f_Z$-almost every $z \in \Z$, it then holds that
	\[
	\bar{K}_1(z, \cdot) = \int f_{X \mid Z}(\df x \mid z) \, f_{Y \mid X}(\df y \mid x) \, f_{X \mid Y}(\df x' \mid y) \, f_{Z \mid X}(\cdot \mid x') = \int_{\Y} f_{Y \mid Z}(\df y \mid z) \, f_{Z \mid Y}(\cdot \mid y).
	\]
	Thus, $\bar{K}_1$ can be seen as a DA algorithm targeting $f_Z$ based on the random element pair $(Z,Y)$.
	One can then apply Lemma~\ref{lem:da-elements} in Appendix~\ref{app:da} to obtain the following.
	
	\begin{proposition} \label{pro:peskun-DA-1}
		Assume that $K$ is the DA algorithm associated with the random element pair $(X,Y)$.
		Suppose that, with respect to some $\sigma$-finite underlying measure on $\Y$, for $f_Z$-almost every $z \in \Z$, $f_{Y \mid Z}(\cdot \mid z)$ has density function $Df_{Y \mid Z}(\cdot \mid z)$.
		Let $\breve{K}_1(z, \df z')$ be an Mtk that is reversible with respect to $f_Z$.
		Suppose that, for $z \in \Z$ and measurable $\mathsf{A} \subset \Z$,
		\[
		\breve{K}_1(z, \mathsf{A} \setminus \{z\}) \leq \int_{\mathsf{A}} \left[ \int_{\Y} \sqrt{Df_{Y \mid Z}(y \mid z) \, Df_{Y \mid Z}(y \mid z')}  \, \df y \right]^2 \, f_Z(\df z').
		\]
		Then, $\bar{K}_1$ satisfies an $(s,b)$-Poincar\'{e} inequality for some $s,b \geq 0$ whenever $\breve{K}_1$ satisfies an $(s,b)$-Poincar\'{e} inequality with the same value of $(s,b)$.
	\end{proposition}

	\subsection{Spectral gap of \texorpdfstring{$\bar{K}_0$}{the overlap kernel} for location and scale mixtures}
	
	We apply Proposition~\ref{pro:peskun-DA-0} to two types of mixture targets.
	These two types respectively encompass \ref{ex:location} and \ref{ex:scale} in the Introduction.
	
	\subsubsection{Location mixture} \label{sssec:Kbar0-location}
	
	Let $\X = \mathbb{R}^d$, and let $m: \Z \to \mathbb{R}^d$ be measurable.
	Let $\xi: \mathbb{R}^d \to \mathbb{R}$ be a differentiable function such that, for some $L \in [0,\infty)$,
	\begin{equation} \label{ine:location-smooth} 
		\xi(x') - \xi(x) - \nabla \xi(x)^{\top} (x'-x) \leq  \frac{L \|x'-x\|^2}{2},  \quad (x,x') \in \mathbb{R}^{d + d}.
	\end{equation}
	Assume that, for $z \in \Z$, $f_{X \mid Z}(\cdot \mid z)$ admits a density function with respect to the Lebesgue measure, given by $x \mapsto Df_{X \mid Z}(x \mid z) = \exp[-\xi(x - m(z))]$.
	That is, $Df_{X \mid Z}(\cdot \mid z)$ is the probability density function $\exp[-\xi(\cdot)]$ shifted by $m(z)$.
	
	The spectral gap of $\bar{K}_0$ heavily depends on the distribution of $m(Z)$ where $Z$ follows the mixing distribution $f_Z$.
	We prove the following proposition in Appendix~\ref{app:Kbar0-location} using Proposition~\ref{pro:peskun-DA-0}.
	
	\begin{proposition} \label{pro:Kbar0-location}
		Consider the context of Section~\ref{sssec:Kbar0-location}.
		Let $\breve{K}_0$ be an Mtk on $\Z$ of the form
		\[
		\breve{K}_0(z, \mathsf{A}) = \int_{\mathsf{A}} \psi(z,z') \, f_Z(\df z') + r(z) \, \ind_{\mathsf{A}}(z),
		\]
		where $r(z) \in \mathbb{R}$ for $z \in \Z$, and $\psi(z,z') = \psi(z',z)$ for $(z,z') \in \Z^2$.
		(Note that $\breve{K}_0$ is reversible with respect to $f_Z$.)
		Assume that $\psi(z,z') \leq \exp[-L\|m(z) - m(z')\|^2/4]$ whenever $z \neq z'$.
		Then $\bar{K}_0$ satisfies an $(s,b)$-Poincar\'{e} inequality whenever $\breve{K}_0$ satisfies an $(s,b)$-Poincar\'{e} inequality with the same value of $(s,b)$.
	\end{proposition}

	An immediate consequence of Proposition~\ref{pro:Kbar0-location} is that, if $z \mapsto m(z)$ is bounded, then $\bar{K}_0$ has a spectral gap.
	
	\begin{corollary} \label{cor:Kbar0-location-bounded}
		Consider the context of Section~\ref{sssec:Kbar0-location}.
		Suppose that there exists $R \geq 0$ such that, for $z, z' \in \Z$, $\|m(z) - m(z')\| \leq R$.
		Then $\mathrm{Gap}(\bar{K}_0) \geq \exp(-LR^2/4)$.
	\end{corollary}
	
	\begin{proof}
		In Proposition~\ref{pro:Kbar0-location}, take $\psi(z,z') = \exp[-L\|m(z)-m(z')\|^2/4]$.
		Then, for $z \in \Z$ and $\mathsf{A} \subset \Z$,
		\[
		\breve{K}_0(z, \mathsf{A} \setminus \{z\}) \geq \exp(-LR^2/4) \, f_Z(\mathsf{A} \setminus \{z\}).
		\]
		Then, for $g \in L_0^2(f_Z)$,
		\[
		\mathcal{E}_{\breve{K}_0}(g) \geq \frac{1}{2} \exp\left(-\frac{LR^2}{4}\right) \int_{\Z^2} f_Z(\df z) \, f_{Z}(\df z') \,[g(z) - g(z')]^2 = \exp\left(-\frac{LR^2}{4}\right) \,\|g\|_{f_Z}^2.
		\]
		Thus, $\mathrm{Gap}(\breve{K}_0) \geq \exp(-LR^2/4)$.
		The desired result then follows from Proposition~\ref{pro:Kbar0-location}.
	\end{proof}
	
	When $\Z$ is finite, one can utilize Proposition~\ref{pro:Kbar0-location} in conjunction with existing techniques and results for Markov chains on discrete state spaces.
	Below is a simple example.
	
	\begin{corollary} \label{cor:Kbar0-location-line}
		Consider the context of Section~\ref{sssec:Kbar0-location}.
		Assume that $\Z = [p]$, where $p\geq2$ is an integer, and that $f_Z(z) = 1/p$ for $z \in [p]$.
		Assume that there exists $\kappa \geq 0$ such that $\|m(z) - m(z+1)\| \leq \kappa$ for $z \in [p-1]$.
		Then
		\[
		\mathrm{Gap}(\bar{K}_0) \geq \frac{2}{p}\left(1-\cos\frac{\pi}{p}\right) \exp\left( - \frac{L \kappa^2}{4} \right).
		\]
	\end{corollary}
	
	\begin{proof}
		In Proposition~\ref{pro:Kbar0-location}, let $\psi(z,z')$ be $\exp(-L\kappa^2/4)$ if $|z-z'| = 1$ and zero otherwise.
		Then $\breve K_0=I-p^{-1}e^{-L\kappa^2/4}\mathsf L_p$, where $\mathsf L_p$ is the Laplacian of the path with $p$ vertices. Its eigenvalues are $2[1-\cos(k\pi/p)]$, $k=0,\ldots,p-1$ \citep[Section~1.4.4]{brouwer2012spectra}. Thus the eigenvalues of $\breve K_0$ are
		\[
		1-\frac{2}{p}e^{-L\kappa^2/4}\left(1-\cos\frac{k\pi}{p}\right),
		\qquad k=0,\ldots,p-1,
		\]
		and its spectral gap is the displayed bound. 
		The desired result then follows from Proposition~\ref{pro:Kbar0-location}.
	\end{proof}
	
	\subsubsection{Scale mixture} \label{sssec:Kbar0-scale}
	
	Let $\X = \mathbb{R}^d$ and $\Z = (0,\infty)$.
	Let $\xi: \mathbb{R}^d \to \mathbb{R}$ be a function such that $x \mapsto \exp[-\xi(x)]$ is a probability density function satisfying $\int_{\mathbb{R}^d} x \, \exp[-\xi(x)] \, \df x = 0$.
	Assume that there are positive numbers $a$ and $L$ such that the eigenvalues of the Hessian matrix of~$\xi$, evaluated at every point in $\mathbb{R}^d$, reside in the interval $[a,L]$.
	Suppose that $f_{X \mid Z}(\cdot \mid z)$ has density function $Df_{X \mid Z}(x \mid z) = z^{-d/2} \exp[-\xi(x/\sqrt{z})]$.
	That is, if a random vector $W$ is distributed according to the probability density function $x \mapsto \exp[-\xi(x)]$, then $\sqrt{z} \,W$ is distributed as $f_{X \mid Z}(\cdot \mid z)$.

	We relate the spectral gap of $\bar{K}_0$ to the tail behavior of $f_Z$, or more precisely, $f_{\log Z}$ defined below.
	Assume that the absolutely continuous part of $f_Z$ with respect to the Lebesgue measure has a density function $Df_Z$.
	With the random variable $Z$ having the distribution $f_Z$, $\log Z$ has distribution $f_{\log Z}(\cdot) = f_Z (\exp(\cdot))$, whose absolutely continuous part has density
	\[
	u \mapsto Df_{\log Z}(u) = e^uDf_Z(e^u) .
	\]
	
	As mentioned in the Introduction, we consider the isoperimetric constant of $f_{\log Z}$, denoted by $\mathrm{Is}(f_{\log Z})$.
	The exact formula is \eqref{eq:is}.
	Based on the isoperimetry of $f_Z$ and conductance-based techniques \citep{lovasz1990mixing,dwivedi2019log,andrieu2024explicit,qin2024convergence}, one can establish the following result, proved in Appendix~\ref{app:Kbar0-scale}.
	
	\begin{proposition} \label{pro:Kbar0-scale}
		In the context of Section~\ref{sssec:Kbar0-scale},
		\[
		\mathrm{Gap}(\bar{K}_0) \geq \frac{1}{128} \left[1 \wedge \frac{\mathrm{Is}(f_{\log Z})^2}{(1+L/a)d} \right].
		\]
	\end{proposition}

	As an example, let $\theta > 0$, and let 
	\begin{equation} \label{eq:ig}
		Df_Z(z) = \frac{(\theta/2)^{\theta/2}}{\Gamma(\theta/2)} z^{-\theta/2-1} \exp \left( - \frac{\theta}{2z} \right), \quad z > 0,
	\end{equation}
	where $\Gamma(\cdot)$ is the gamma function.
	This is the inverse-gamma distribution with shape parameter $\theta/2$ and scale parameter $\theta/2$.
	A small value of $\theta$ results in $Df_Z$ having a heavy right tail, which in turn results in $Df_X$ having heavy tails. 
	We prove the following proposition in Appendix~\ref{app:Is-ig}.

	\begin{proposition} \label{pro:Is-ig}
		Suppose that $f_Z$ has density function~\eqref{eq:ig}, where $\theta>0$.
		Then
		\[
		\mathrm{Is}(f_{\log Z})
		\geq\frac{(\theta/2)^{\theta/2}\exp(-\theta/2)}{\Gamma(\theta/2)}
		\geq\frac{\theta}{4\sqrt{\theta+2}}.
		\]
		Thus, by Proposition~\ref{pro:Kbar0-scale},
		\[
		\mathrm{Gap}(\bar K_0)\geq\frac1{128}
		\left[1\wedge\frac{\theta^2}{16(\theta+2)(1+L/a)d}\right].
		\]
	\end{proposition}

	\section{Applications} \label{sec:applications}
	
	\subsection{Overview}
	
	We apply the decomposition results in Section \ref{sec:decomposition} to Gaussian random-walk Metropolis--Hastings, the proximal sampler, and random-scan Gibbs sampling. 
	Existing results for well-conditioned targets give spectral gap estimates for $K_z$, $z \in \Z$. 
	We combine these local results with bounds involving $\bar{K}_0$ or $\bar{K}_1$ to derive Poincar\'{e} inequalities for~$K$.
	
	\subsection{Random-walk Metropolis--Hastings Algorithm} \label{ssec:MH}
	
	\subsubsection{Overview} \label{sssec:MH-overview}
	
	Let $\X=\mathbb R^d$ and let $f_X$ have Lebesgue density $Df_X$. Consider \ref{ex:MH} with proposal density
	\[
	q(x,x')=(2\pi\sigma^2)^{-d/2}\exp[-\|x'-x\|^2/(2\sigma^2)],\qquad \sigma>0,
	\]
	and acceptance probability $1\wedge [Df_X(x')/Df_X(x)]$. For each component $z \in \Z$ associated with the conditional density $Df_{X\mid Z}(\cdot\mid z)$, let $K_z$ correspond to the MH algorithm targeting $f_{X \mid Z}(\cdot \mid z)$ with the same proposal $q(x,x')$ and acceptance probability $1\wedge Df_{X\mid Z}(x'\mid z)/Df_{X\mid Z}(x\mid z)$. Proposition~\ref{pro:condition-MH} verifies Condition~\ref{cond:A}.
	
	We use the following consequence of Theorem~1 of \cite{andrieu2024explicit}. 
	
	\begin{proposition} \label{pro:MH-concave}
		Suppose $-\log Df_{X\mid Z}(\cdot\mid z)$ is differentiable, $a_z$-strongly convex and $L_z$-smooth, with $0<a_z\leq L_z$. Then
		\[
		\mathrm{Gap}(K_z)\geq 10^{-4}L_zd\sigma^2\exp(-2L_zd\sigma^2)\frac{a_z}{L_zd}.
		\]
	\end{proposition}
	
	\begin{proof}
		Apply the lower bound in Theorem~1 of \cite{andrieu2024explicit} with strong-convexity constant $a_z$ and smoothness constant $L_z$, and replace its numerical constant $1.972\times10^{-4}$ by $10^{-4}$. As noted immediately after that theorem, twice continuous differentiability is required only for its upper bound, so the differentiability assumed here suffices for the lower bound.
	\end{proof}
	
	\subsubsection{Location mixture} \label{sssec:MH-location}
	
	Let $Df_{X\mid Z}(x\mid z)=\exp[-\xi(x-m(z))]$, where $m:\Z\to\mathbb R^d$ is measurable, $\exp(-\xi)$ is a probability density, and $\xi$ is differentiable with
	\begin{equation} \label{ine:location-convex-smooth}
		\frac{a\|x'-x\|^2}{2}\leq \xi(x')-\xi(x)-\nabla\xi(x)^\top(x'-x)\leq\frac{L\|x'-x\|^2}{2},\qquad 0<a\leq L.
	\end{equation}
	
	\begin{proposition} \label{pro:Kz-MH-location}
		If $\sigma=\varsigma/\sqrt{Ld}$, where $\varsigma>0$, then every component kernel satisfies
		\[
		\mathrm{Gap}(K_z)\geq10^{-4}\varsigma^2e^{-2\varsigma^2}\frac{a}{Ld}.
		\]
	\end{proposition}
	
	\begin{proof}
		The potential of $f_{X \mid Z}(\cdot \mid z)$ is $-\log Df_{X \mid Z}(x \mid z) = \xi(x-m(z))$.
		Translation preserves the strong-convexity and smoothness bounds in~\eqref{ine:location-convex-smooth}. Apply Proposition~\ref{pro:MH-concave} to $K_z$ with $a_z=a$, $L_z=L$, and $Ld\sigma^2=\varsigma^2$.
		This gives the desired result.
	\end{proof}
	
	We may now state a spectral gap bound for~$K$ targeting the location mixture.
	
	\begin{proposition} \label{pro:MH-location}
		Assume $\sigma = \varsigma/\sqrt{Ld}$, where $\varsigma > 0$.
		Let $\breve K_0$ be an $f_Z$-reversible Markov kernel of the form
		\[
		\breve K_0(z,\mathsf A)=\int_{\mathsf A}\psi(z,z')f_Z(\df z')+r(z)\ind_{\mathsf A}(z),
		\]
		where $\psi$ is symmetric and $\psi(z,z')\leq\exp[-L\|m(z)-m(z')\|^2/4]$ for $z\ne z'$. Then
		\begin{equation} \nonumber
			\mathrm{Gap}(K)\geq10^{-4}\varsigma^2e^{-2\varsigma^2}\frac{a}{Ld}\,\mathrm{Gap}(\breve K_0).
		\end{equation}
		In particular:
		\begin{enumerate}
			\item[(i)] If $\sup_{z,z'}\|m(z)-m(z')\|\leq R$, then $\mathrm{Gap}(K)\geq10^{-4}\varsigma^2e^{-2\varsigma^2}\frac{a}{Ld}e^{-LR^2/4}$.
			\item[(ii)] If $\Z=[p]$, $p\geq2$, $f_Z(z)=1/p$, and $\|m(z)-m(z+1)\|\leq\kappa$, then
			\[
			\mathrm{Gap}(K)\geq10^{-4}\varsigma^2e^{-2\varsigma^2}\frac{a}{Ld}\frac{2}{p}\bigl[1-\cos(\pi/p)\bigr]e^{-L\kappa^2/4}.
			\]
		\end{enumerate}
	\end{proposition}
	
	\begin{proof}
		Combine Corollary~\ref{cor:uniform-gap}, Proposition~\ref{pro:Kz-MH-location}, and Proposition~\ref{pro:Kbar0-location}. The special cases follow from Corollaries~\ref{cor:Kbar0-location-bounded} and~\ref{cor:Kbar0-location-line}.
	\end{proof}

	Spectral gap estimates can be translated to mixing-time bounds using Corollary \ref{cor:mixing-sb}.
	
	\begin{corollary} \label{cor:mixing-MH-location}
		Fix $\varsigma>0$.
		Let $\epsilon > 0$, $\mu\in L_*^2(f_X)$.
		Under the respective assumptions of Proposition~\ref{pro:MH-location},
		\[
		t(\epsilon;\mu,K)=
		\begin{cases}
			O\bigl( (L/a)d\,e^{LR^2/4} \log(\|\mu-f_X\|_{f_X}/\epsilon) \bigr),&\text{in (i)},\\
			O\bigl( (L/a)d\,p^3e^{L\kappa^2/4} \log(\|\mu-f_X\|_{f_X}/\epsilon) \bigr),&\text{in (ii)}.
		\end{cases}
		\]
		In (i), the big O notation depicts asymptotics when $\|\mu - f_X\|_{f_X}/\epsilon$, $L/a$, $d$, and $LR^2$ tend to infinity; in (ii), the big O notation depicts asymptotics when $\|\mu - f_X\|_{f_X}/\epsilon$, $L/a$, $d$, $p$, $L\kappa^2$ tend to infinity.
	\end{corollary}
	
	\begin{proof}
		The Gaussian-proposal Metropolis kernel is positive semi-definite,
		as recalled in Section~\ref{sec:Poincare}.
		Apply Corollary~\ref{cor:mixing-sb} with $b=0$ and $s$ equal to the
		reciprocal of the corresponding gap lower bound in
		Proposition~\ref{pro:MH-location}.
		%Use $-\log(1-1/s)\geq1/s$ and
		%$\log(\|\mu-f_X\|_{f_X}/(2\epsilon))
		%=O(\log\|\mu-f_X\|_{f_X})$ in the stated limit.
		In case~(ii), $1-\cos(\pi/p) = 2 \sin^2[\pi/(2p)] \geq2/p^2$ for $p\geq 1$, giving the
		factor $p^3$.
	\end{proof}
	
	%The $p^3$ factor in this comparison comes from the adjacent transition probability $e^{-L\kappa^2/4}/p$ and the $p^{-2}$ spectral gap of a path. It is an upper-bound consequence of this particular comparison and is not an assertion that the original sampler necessarily requires cubic time in $p$.
	
	\subsubsection{Scale mixture} \label{sssec:MH-scale}
	
	Consider the scale mixture scenario in Section \ref{sssec:Kbar0-scale}.
	That is, let $\Z=(0,\infty)$ and
	\[
	Df_{X\mid Z}(x\mid z)=z^{-d/2}\exp[-\xi(x/\sqrt z)],
	\]
	where $\exp(-\xi)$ is a centered probability density and $\xi$ is twice continuously differentiable with $aI_d\preceq\nabla^2\xi\preceq LI_d$, $0<a\leq L$. 
	
	\begin{proposition} \label{pro:Kz-MH-scale}
		For $\sigma=\varsigma/\sqrt{Ld}$, where $\varsigma>0$,
		\[
		\mathrm{Gap}(K_z)\geq10^{-4}\frac{\varsigma^2}{z}e^{-2\varsigma^2/z}\frac{a}{Ld}.
		\]
	\end{proposition}
	
	\begin{proof}
		The potential of $f_{X \mid Z}(\cdot \mid z)$ is $x \mapsto \xi(x/\sqrt z)+(d/2)\log z$.
		It is $(a/z)$-strongly convex and $(L/z)$-smooth. Substitute $a_z=a/z$, $L_z=L/z$, and $L_zd\sigma^2=\varsigma^2/z$ into Proposition~\ref{pro:MH-concave}.
	\end{proof}
	
	We may now state a Poincar\'{e} inequality for $K$ targeting the scale mixture.
	
	\begin{proposition} \label{pro:MH-scale}
		Assume that $\mathrm{Is}(f_{\log Z}) > 0$.
		Under the assumptions of Proposition \ref{pro:Kz-MH-scale}, $K$ satisfies a weak Poincar\'{e} inequality with rate
		\[
		\begin{aligned}
			\beta(s)=\frac14\wedge\biggl\{&32\left[1\vee\frac{(1+L/a)d}{\mathrm{Is}(f_{\log Z})^2}\right] \Pr\left(128\left[1\vee\frac{(1+L/a)d}{\mathrm{Is}(f_{\log Z})^2}\right]
			\frac{10^4 (L/a) d Z}{\varsigma^2}e^{2\varsigma^2/Z}>s\right)\biggr\},\qquad s>0,
		\end{aligned}
		\]
		where $Z\sim f_Z$.
	\end{proposition}
	
	\begin{proof}
		Apply Part~(i) of Proposition~\ref{pro:main}, using Propositions~\ref{pro:Kz-MH-scale} and~\ref{pro:Kbar0-scale}. The trivial variance bound $\operatorname{Var}(g)\leq\|g\|_{\rm osc}^2/4$ permits the minimum with $1/4$.
	\end{proof}
	
	\begin{corollary} \label{cor:MH-scale}
		Assume the conditions of Proposition \ref{pro:Kz-MH-scale}, with $\varsigma > 0$ fixed.
		Suppose further that $\mathbb{E}(Z e^{2\varsigma^2/Z}) < \infty$, and that $0 < \mathrm{Is}(f_{\log Z})^2 \leq (1+L/a) d$.
		Then $K$ satisfies a weak Poincar\'{e} inequality with rate
		\[
		\beta(s) = \frac{1}{4} \wedge \left[ \frac{ 4.1 \times 10^7 (1+L/a)^2 (L/a) d^3}{\mathrm{Is}(f_{\log Z})^4} \, \mathbb{E}\left( \frac{Z}{\varsigma^2} e^{2\varsigma^2/Z} \right) s^{-1} \right].
		\]
		Consequently, for $\epsilon > 0$ and $\mu \in L_*^2(f_X)$,
		as $L/a$, $d$, $\|\mu - f_X \|_{f_X} /\epsilon \to \infty$,
		\[
		t(\epsilon; \mu, K) = O \left( \left( \frac{L}{a} \right)^3 d^3 \frac{\|\mu-f_X\|_{f_X}^2}{\epsilon^2} \right).
		\]
	\end{corollary}
	\begin{proof}
		The weak Poincar\'{e} inequality is obtained by combining Proposition \ref{pro:MH-scale} with Markov inequality.
		The mixing-time bound follows from Corollary \ref{cor:mixing-polynomial}.
	\end{proof}
	
	%\begingroup\color{blue}
	%\begin{proposition} \label{pro:MH-ig}
	%    Let $f_Z$ be the inverse-gamma law~\eqref{eq:ig}, with $\theta>0$.
	%    Under the assumptions of Proposition \ref{pro:Kz-MH-scale}, a valid weak Poincar\'{e} rate, for $0<\ell<\min\{\theta/2,\theta/(4\varsigma^2)\}$, is
	%    \[
	%    \begin{aligned}
		%    \beta(s)=\frac14\wedge\biggl\{&32\cdot128^\ell
		%    \left[1\vee\frac{16(\theta+2)(1+L/a)d}{\theta^2}\right]^{\ell+1}
		%    \left(\frac{10^4Ld}{a\varsigma^2s}\right)^\ell\\
		%    &\times\frac{(\theta/2)^{\theta/2}\Gamma(\theta/2-\ell)}
		%    {\Gamma(\theta/2)(\theta/2-2\ell\varsigma^2)^{\theta/2-\ell}}\biggr\}.
		%    \end{aligned}
	%    \]
	%\end{proposition}
	%
	%\begin{proof}
	%In Proposition~\ref{pro:MH-scale}, Markov's inequality permits the simpler rate
	%\[
	%\beta(s)=\frac14\wedge\left\{
	%32\cdot128^\ell\left[1\vee\frac{(1+L/a)d}{\mathrm{Is}(f_{\log Z})^2}\right]^{\ell+1}
	%\left(\frac{10^4Ld}{a\varsigma^2s}\right)^\ell
	%\mathrm E\bigl(Z^\ell e^{2\ell\varsigma^2/Z}\bigr) \right\}.
	%\]
	%Using the inverse-gamma density~\eqref{eq:ig} gives
	%\[
	%\begin{aligned}
	%\mathrm E(Z^\ell e^{2\ell\varsigma^2/Z})
	%&=\frac{(\theta/2)^{\theta/2}}{\Gamma(\theta/2)}
	%\int_0^\infty z^{\ell-\theta/2-1}
	%e^{-(\theta/2-2\ell\varsigma^2)/z}\,\df z\\
	%&=\frac{(\theta/2)^{\theta/2}\Gamma(\theta/2-\ell)}
	%{\Gamma(\theta/2)(\theta/2-2\ell\varsigma^2)^{\theta/2-\ell}}.
	%\end{aligned}
	%\]
	%The restrictions on $\ell$ ensure that the integral is finite. Finally, Proposition~\ref{pro:Is-ig} gives $\mathrm{Is}(f_{\log Z})^{-2}\leq16(\theta+2)/\theta^2$.
	%\end{proof}
	%\endgroup

	Recall that, when $f_Z$ is the inverse gamma law~\eqref{eq:ig} with $\theta > 0$, a smaller $\theta$ result in heavier tails.
	The following corollary, which is proved in Appendix \ref{app:mixing-MH-scale}, establishes a relationship between the mixing time and~$\theta$.
	Corollary \ref{cor:mixing-MH-scale} is not a direct corollary of Corollary \ref{cor:MH-scale} as it requires a more intricate analysis specific to inverse gamma distributions.

	\begin{corollary} \label{cor:mixing-MH-scale}
		Let $f_Z$ be the inverse-gamma law~\eqref{eq:ig}, with $\theta>0$, and assume the conditions of
		Proposition \ref{pro:Kz-MH-scale}, with $\varsigma \in (0,1/\sqrt{2})$ fixed.
		Let $\mu\in L_*^2(f_X)$ and $\epsilon>0$. As
		$L/a,d,\|\mu-f_X\|_{f_X}/\epsilon \to\infty$ and $\theta\downarrow0$,
		\[
		t(\epsilon;\mu,K)=O\!\left(
		\frac{(L/a)^2 d^2}{\theta^2}
		\left[\frac{1024(1+L/a)d}
		{\theta^2} \frac{\|\mu-f_X\|_{f_X}^2}{\epsilon^2}\right]^{2/\theta}
		\right).
		\]
	\end{corollary}

	%For fixed $\theta>0$, the same argument gives the endpoint weak Poincar\'{e} rate $\beta(s)=O(s^{-\theta/2})$ and, by Proposition~\ref{pro:polynomial},
	%$\mathrm{TV}(\mu K^t,f_X)=O(t^{-\theta/4})$ as $t\to\infty$.
	%The tail coefficient in~\eqref{ine:MH-ig-tail} is finite and continuous at $\theta=2$; the first-moment singularity at that value is therefore absent from this bound.

	\subsection{Proximal Sampler} \label{ssec:proximal}
	
	\subsubsection{Overview} \label{sssec:proximal-overview}
	
	Let $\X=\Y=\mathbb R^d$ and $f_{Y\mid X}(\cdot\mid x)=\mathrm N_d(x,\eta I_d)$, $\eta>0$. The proximal sampler alternates this Gaussian draw with an exact draw from
	\[
	Df_{X\mid Y}(x\mid y)\propto Df_X(x)\exp[-\|y-x\|^2/(2\eta)].
	\]
	Its kernel is $K(x,\df x')=\int f_{Y\mid X}(\df y\mid x)f_{X\mid Y}(\df x'\mid y)$. Following Proposition~\ref{pro:condition-DA}, define $K_z$ by replacing the second conditional with $f_{X\mid Y,Z}(\cdot\mid y,z)$. Then Condition~\ref{cond:A} holds and $K_z$ is the proximal sampler for the $z$th component.
	
	\begin{proposition} \label{pro:Kz-concave} \citep{chen2022improved}
		If $-\log Df_{X\mid Z}(\cdot\mid z)$ is $a_z$-strongly convex, $a_z>0$, then
		\[
		\mathrm{Gap}(K_z)\geq\frac{a_z\eta}{1+a_z\eta}.
		\]
	\end{proposition}
	
	\begin{proof}
		Strong log-concavity gives a differential Poincar\'{e} inequality for the component distribution $f_{X \mid Z}(\cdot \mid z)$ with constant $1/a_z$. Theorem~4 of \cite{chen2022improved} therefore contracts the $\chi^2$ divergence by a factor $(1+a_z\eta)^{-2}$ in one iteration. Apply this to initial densities $1+\varepsilon g$, for bounded mean-zero $g$ and sufficiently small $\varepsilon>0$, to obtain
		\[
		\|K_zg\|_{f_{X\mid Z}(\cdot\mid z)}
		\leq\frac{1}{1+a_z\eta}\|g\|_{f_{X\mid Z}(\cdot\mid z)}.
		\]
		Bounded mean-zero functions are dense in the space of $L^2$ mean-zero functions, so the same inequality holds for every mean-zero $g$ in $L^2(f_{X \mid Z}(\cdot \mid z))$. Reversibility then gives $\mathrm{Gap}(K_z)\geq1-(1+a_z\eta)^{-1}$.
	\end{proof}
	
	Write $K(\eta)$ to emphasize~$K$ depends on~$\eta$. The next result is proved in Appendix~\ref{app:ordering-eta}.
	
	\begin{proposition} \label{pro:ordering-eta}
		If $0<\eta\leq\eta'$, then $\mathcal E_{K(\eta)}(g)\leq\mathcal E_{K(\eta')}(g)$ for every $g\in L^2(f_X)$. Consequently, every $(s,b)$-Poincar\'{e} inequality for $K(\eta)$ holds for $K(\eta')$ with the same constants.
	\end{proposition}
	
	A larger $\eta$ improves the mixing of the proximal sampler in this sense, but may make the conditional draw from $f_{X \mid Y}$ harder. For fixed $y$, $f_{X \mid Y}(\cdot \mid y)$ approaches the target distribution $f_X(\cdot)$ as $\eta\to\infty$, making implementing the proximal sampler as difficult as the original sampling problem.
	Section~\ref{sssec:proximal-metropolis} studies the mixing properties of a Metropolis implementation.
	
	\subsubsection{Location mixture} \label{sssec:proximal-location}
	
	Assume the location-mixture conditions in Section~\ref{sssec:MH-location}.
	
	By Proposition \ref{pro:Kz-concave}, the following holds.
	
	\begin{proposition} \label{pro:Kz-proximal-location}
		For every $z \in \Z$, $\mathrm{Gap}(K_z)\geq a\eta/(1+a\eta)$.
	\end{proposition}
	
	\begin{proof}
		The potential of $f_{X \mid Z}(\cdot \mid z)$ is $x\mapsto\xi(x-m(z))$, whose strong-convexity constant is $a$. The desired result then follows from Proposition~\ref{pro:Kz-concave} with $a_z=a$.
	\end{proof}
	
	We will combine Proposition~\ref{pro:Kz-proximal-location} with a Poincar\'{e} inequality for $\bar{K}_1$.
	To this end, we state the following result, which is proved in Appendix~\ref{app:Kbar1-proximal-location}.
	
	\begin{proposition} \label{pro:Kbar1-proximal-location}
		Let $\breve K_1$ be an $f_Z$-reversible Markov kernel of the form
		\[
		\breve K_1(z,\mathsf A)=\int_{\mathsf A}\psi(z,z')f_Z(\df z')+r(z)\ind_{\mathsf A}(z),
		\]
		where $\psi$ is symmetric and $\psi(z,z')\leq\exp[-L\|m(z)-m(z')\|^2/[4(1+L\eta)]]$ for $z\ne z'$. Every $(s,b)$-Poincar\'{e} inequality for $\breve K_1$ also holds for $\bar K_1$.
	\end{proposition}

	Proposition \ref{pro:Kbar1-proximal-location} is applied to two concrete scenarios.
	(i) The $m(z)$'s lie within a ball of some diameter $R$, and (ii) $\Z = [p]$ with $f_Z(z) = 1/p$, and $\|m(z) - m(z+1)\| \leq \kappa$, where $\kappa \geq 0$. 
	
	\begin{proposition} \label{pro:Kbar1-proximal-location-1}
		The following lower bounds hold:
		\begin{enumerate}
			\item[(i)] If $\sup_{z,z'}\|m(z)-m(z')\|\leq R$, then $\mathrm{Gap}(\bar K_1)\geq e^{-LR^2/[4(1+L\eta)]}$.
			\item[(ii)] If $\Z=[p]$, $p\geq2$, $f_Z(z)=1/p$ and $\|m(z)-m(z+1)\|\leq\kappa$, then
			\[
			\mathrm{Gap}(\bar K_1)\geq\frac2p[1-\cos(\pi/p)]e^{-L\kappa^2/[4(1+L\eta)]}.
			\]
		\end{enumerate}
	\end{proposition}
	
	\begin{proof}
		For part~(i), choose $\psi(z,z')=e^{-LR^2/[4(1+L\eta)]}$ in Proposition~\ref{pro:Kbar1-proximal-location}. The resulting kernel refreshes from $f_Z$ with this probability, so its Dirichlet form is this probability times the variance under $f_Z$. For part~(ii), choose
		\[
		\psi(z,z')=e^{-L\kappa^2/[4(1+L\eta)]}\ind_{\{|z-z'|=1\}}.
		\]
		This is the path kernel in the proof of Corollary~\ref{cor:Kbar0-location-line}, with $e^{-L\kappa^2/4}$ replaced by $e^{-L\kappa^2/[4(1+L\eta)]}$, so its spectral gap is the displayed bound. Proposition~\ref{pro:Kbar1-proximal-location} transfers both bounds to $\bar K_1$.
	\end{proof}

	\begin{proposition} \label{pro:proximal-location}
		For $\breve K_1$ given in Proposition~\ref{pro:Kbar1-proximal-location},
		\begin{equation} \nonumber
			\mathrm{Gap}(K)\geq \frac{a\eta}{1+a\eta}\,\mathrm{Gap}(\breve K_1).
		\end{equation}
		In particular, under the respective conditions of Proposition~\ref{pro:Kbar1-proximal-location-1},
		\[
		\mathrm{Gap}(K)\geq
		\begin{cases}
			\dfrac{a\eta}{1+a\eta}e^{-LR^2/[4(1+L\eta)]},&\text{in (i)},\\
			\dfrac{a\eta}{1+a\eta}\dfrac2p[1-\cos(\pi/p)]e^{-L\kappa^2/[4(1+L\eta)]},&\text{in (ii)}.
		\end{cases}
		\]
	\end{proposition}
	
	\begin{proof}
		Apply the final bound in Corollary~\ref{cor:uniform-gap}, which extends \citet[Proposition~2.4]{CLY2026}, together with Propositions~\ref{pro:Kz-proximal-location} to~\ref{pro:Kbar1-proximal-location-1}.
	\end{proof}
	
	%\begin{remark}
	%The independent bound from $\bar K_0$ is also available: in the bounded-diameter case,
	%\[
	%\mathrm{Gap}(K)\geq\max\left\{\frac{a\eta}{1+a\eta}e^{-LR^2/4},
	%\frac{a\eta}{1+a\eta}e^{-LR^2/[4(1+L\eta)]}\right\}.
	%\]
	%The second bound dominates the first after applying the DA product bound.
	%\end{remark}

	\begin{remark}
		\citet[Theorem~2.5]{CLY2026} gives the bounded-diameter estimate $[a\eta/(1+a\eta)]e^{-R^2/(4\eta)}$ for $\mathrm{Gap}(K)$. 
		Under our additional $L$-smoothness assumption, Proposition~\ref{pro:proximal-location}(i) improves its overlap exponent, since $L/(1+L\eta)<1/\eta$.
	\end{remark}
	
	\begin{corollary} \label{cor:mixing-proximal-location}
		Let $\epsilon>0$ and let $\mu\in L_*^2(f_X)$.
		Under the respective conditions of Proposition~\ref{pro:Kbar1-proximal-location-1},
		\[
		t(\epsilon;\mu,K)=
		\begin{cases}
			O\bigl([1+(a\eta)^{-1}]e^{LR^2/[4(1+L\eta)]} \log(\|\mu-f_X\|_{f_X}/\epsilon) \bigr),&\text{in (i)},\\
			O\bigl([1+(a\eta)^{-1}]p^3e^{L\kappa^2/[4(1+L\eta)]} \log (\|\mu-f_X\|_{f_X}/\epsilon) \bigr),&\text{in (ii)}.
		\end{cases}
		\]
		In (i), the big $O$ notation depicts asymptotics when $\|\mu - f_X\|_{f_X}/\epsilon$, $(a\eta)^{-1}$, and $LR^2/(1+L\eta)$ tend to infinity;
		in (ii), the big $O$ notation depicts asymptotics when $\|\mu - f_X\|_{f_X}/\epsilon$, $(a\eta)^{-1}$, $p$, and $L\kappa^2/(1+L\eta)$ tend to infinity.
	\end{corollary}
	
	\begin{proof}
		The proximal kernel is positive semi-definite.
		Apply Corollary~\ref{cor:mixing-sb} with $b=0$ and $s$ equal to the
		reciprocal of the relevant gap lower bound in
		Proposition~\ref{pro:proximal-location}.
	\end{proof}
	
	\subsubsection{Scale mixture} \label{sssec:proximal-scale}
	
	Assume the centered scale-mixture conditions of Section~\ref{sssec:Kbar0-scale}.
	
	\begin{proposition} \label{pro:Kz-proximal-scale}
		For every $z>0$,
		\[
		\mathrm{Gap}(K_z)\geq\frac{a\eta}{z+a\eta}.
		\]
	\end{proposition}
	
	\begin{proof}
		The potential of $f_{X \mid Z}(\cdot \mid z)$ is $x\mapsto\xi(x/\sqrt z)+(d/2)\log z$, whose Hessian is bounded below by $(a/z)I_d$. Apply Proposition~\ref{pro:Kz-concave} with $a_z=a/z$.
	\end{proof}
	
	\begin{proposition} \label{pro:proximal-scale}
		Assume that $\mathrm{Is}(f_{\log Z}) > 0$.
		Then $K$ satisfies a weak Poincar\'{e} inequality with
		\[
		\begin{aligned}
			\beta(s)=\frac14\wedge\biggl\{&32\left[1\vee\frac{(1+L/a)d}{\mathrm{Is}(f_{\log Z})^2}\right] \Pr\left(128\left[1\vee\frac{(1+L/a)d}{\mathrm{Is}(f_{\log Z})^2}\right]
			\frac{Z+a\eta}{a\eta}>s\right)\biggr\},\qquad s>0.
		\end{aligned}
		\]
	\end{proposition}
	
	\begin{proof}
	Apply the same argument as in Proposition~\ref{pro:MH-scale}, and use Proposition \ref{pro:Kz-proximal-scale}.
	\end{proof}
	%
	%\begin{remark} \label{rem:proximal-scale}
	%In Proposition~\ref{pro:proximal-scale}, for any $\ell>0$ with $\mathrm E[(Z+a\eta)^\ell]<\infty$, a valid simpler rate is
	%\[
	%\beta(s)=\frac14\wedge\left\{
	%\frac{32\cdot128^\ell}{(a\eta s)^\ell}
	%\left[1\vee\frac{(1+L/a)d}{\mathrm{Is}(f_{\log Z})^2}\right]^{\ell+1}
	%\mathrm E[(Z+a\eta)^\ell]\right\}.
	%\]
	%This follows by applying Markov's inequality to $(Z+a\eta)^\ell$ in Proposition~\ref{pro:proximal-scale}. For the inverse-gamma law~\eqref{eq:ig}, the density decays as $z^{-\theta/2-1}$ at infinity and exponentially in $1/z$ at zero, so these moments are finite for $0<\ell<\theta/2$.
	%\end{remark}
	
	\begin{corollary}
		Assume that $\mathbb{E}(Z) < \infty$, and that $0 < \mathrm{Is}(f_{\log Z})^2 \leq (1+L/a)d$.
		Then $K$ satisfies a weak Poincar\'{e} inequality with
		\[
		\beta(s) = \frac{1}{4} \wedge \left[ \frac{4096 (1+L/a)^2 d^2}{\mathrm{Is}(f_{\log Z})^4} \, \mathbb{E}\left( \frac{Z}{a\eta} + 1 \right) s^{-1} \right].
		\]
		Consequently, for $\epsilon > 0$ and $\mu \in L_*^2(f_X)$, as $L/a$, $d$, $(a\eta)^{-1} \to \infty$,
		\[
		t(\epsilon; \mu, K) = O \left( \left( \frac{L}{a} \right)^2 d^2 (a\eta)^{-1} \frac{\|\mu-f_X\|_{f_X}^2}{\epsilon^2} \right).
		\]
	\end{corollary}
	\begin{proof}
		The weak Poincar\'{e} inequality follows from Proposition \ref{pro:proximal-scale} and Markov's inequality.
		The mixing-time bound follows from Corollary \ref{cor:mixing-polynomial}.
	\end{proof}
	
	For illustration, consider the case where $f_Z$ is an inverse gamma distribution.
	The following corollary is proved in Appendix \ref{app:mixing-proximal-scale}.

	\begin{corollary} \label{cor:mixing-proximal-scale}
		Let $f_Z$ have the inverse-gamma law~\eqref{eq:ig} with $\theta>0$.
		Let $\epsilon > 0$, and $\mu\in L_*^2(f_X)$.
		As
		$L/a,d,\|\mu-f_X\|_{f_X}/\epsilon \to\infty$, $a\eta\to0$, and $\theta\downarrow0$,
		\[
		t(\epsilon;\mu,K)=O\!\left(
		\frac{(L/a)d}{a\eta \, \theta}
		\left[\frac{1024(1+L/a) \, d  }
		{\theta^2} \frac{\|\mu-f_X\|_{f_X}^2}{\epsilon^2}\right]^{2/\theta}
		\right).
		\]
	\end{corollary}

	\subsubsection{Randomized-step MALA implementation} \label{sssec:proximal-metropolis}
	
	In practice, the proximal sampler can be difficult to implement because exact draws from $f_{X\mid Y}(\cdot\mid y)$ are expensive or infeasible.
	A hybrid proximal sampler can be applied in this case.
	
	Let $Q_y$ be reversible with respect to $f_{X\mid Y}(\cdot\mid y)$, and replace the exact conditional draw by one $Q_y$ step. The resulting kernel
	\[
	H(x,\df x')=\int f_{Y\mid X}(\df y\mid x)Q_y(x,\df x')
	\]
	is $f_X$-reversible and is positive semi-definite whenever each $Q_y$ is. The following comparison is the weak Poincar\'{e} form of the hybrid data-augmentation argument; see also \cite{qin2025spectral}.
	
	\begin{proposition} \label{pro:hybrid}
		Suppose that, for $f_Y$-almost every $y$, $Q_y$ satisfies an $(s,b_y)$-Poincar\'{e} inequality, where $s\geq0$, $b_y\geq0$ and $y\mapsto b_y$ is measurable. If $K$ satisfies an $(s',b')$-Poincar\'{e} inequality, then $H$ satisfies one with
		\[
		s''=ss',\qquad b''=b'+s'\int b_y f_Y(\df y),
		\]
		assuming the integral is finite.
		In particular, $\mathrm{Gap}(H)\geq \delta\,\mathrm{Gap}(K)$ whenever $\mathrm{Gap}(Q_y)\geq \delta>0$ almost surely.
	\end{proposition}
	
	\begin{proof}
		Let $g\in L^2(f_X)$. Applying the Poincar\'{e} inequality for $Q_y$ to $g_y(x) = g(x)-\mathrm E[g(X)\mid Y=y]$ yields
		\[
		\operatorname{Var}(g(X) \mid Y = y) \leq \frac{s}{2} \int_{\X^2} f_{X \mid Y}(\df x \mid y) \, Q_y(x, \df x') [g(x) - g(x')]^2 + b_y \|g_y\|_{\mathrm{osc}}^2.
		\]
		It can be checked that $\mathcal E_K(g)=\mathrm E[\operatorname{Var}(g(X)\mid Y)]$, and, for $f_Y$-a.e. $y \in \Y$, $\|g_y\|_{\mathrm{osc}} \leq \|g\|_{\mathrm{osc}}$.
		In the most recent display, integrating with respect to $f_Y$ then gives
		\[
		\mathcal E_K(g)\leq s\mathcal E_H(g)+\left(\int b_y f_Y(\df y)\right)\|g\|_{\rm osc}^2.
		\]
		Insert this into the assumed Poincar\'{e} inequality for $K$ to obtain the desired result.
	\end{proof}
	
	Throughout Section~\ref{sssec:proximal-metropolis}, take $Q_y$ to be the lazy randomized-step Metropolis-adjusted Langevin algorithm (MALA) analyzed in \cite{qin2026randomizedmala}, targeting the conditional potential $-\log Df_X(x)+\|x-y\|^2/(2\eta)$. With probability $1/2$, retain $x$. Otherwise, draw $h$ uniformly from $(0,h_{\max})$ and propose
	\[
	x'\sim\mathrm N_d\!\left(x+h\nabla\log Df_X(x)-\frac{h}{\eta}(x-y),\,2hI_d\right).
	\]
	Apply the Metropolis--Hastings acceptance rule for $f_{X\mid Y}(\cdot\mid y)$ using the Gaussian proposal density at the sampled $h$. A fresh $h$ is drawn at each such update. The holding probability ensures that $Q_y$ is positive semi-definite.
	
	We now give a Poincar\'{e} inequality involving $Q_y$.
	
	\begin{proposition} \label{pro:proximal-semi}
		Suppose that $-\log Df_X$ is differentiable and, for some real $a'\leq L'$,
		\begin{equation} \label{ine:semi}
			\frac{a'\|x'-x\|^2}{2}\leq \log Df_X(x)-\log Df_X(x')+\nabla\log Df_X(x)^\top(x'-x)
			\leq\frac{L'\|x'-x\|^2}{2}, \quad x,x' \in \mathbb{R}^d.
		\end{equation}
		Choose $\eta>0$ so that $\eta^{-1} \geq L'-2a'$, and set
		\[
		h_{\max}=\frac{1}{(L'+\eta^{-1})\sqrt{d(1+\log d)}}.
		\]
		Then, for a universal constant $c\in(0,\infty)$, uniformly in $y$,
		\[
		\mathrm{Gap}(Q_y)\geq\frac{c}{\sqrt{d(1+\log d)}}.
		\]
	\end{proposition}
	
	\begin{proof}
		For $y \in \mathbb{R}^d$, the potential of the conditional distribution $Df_{X \mid Y}(x \mid y)$ is 
		\[
		-\log Df_{X \mid Y}(x \mid y) = - \log Df_X(x) + \frac{\|x-y\|^2}{2\eta} + \text{term independent of }x.
		\]
		By \eqref{ine:semi}, for $x, x' \in \mathbb{R}^d$,
		\[
		\begin{aligned}
			\frac{a'+\eta^{-1}}{2} \|x-x'\|^2 \leq& \log Df_{X \mid Y}(x \mid y) - \log Df_{X \mid Y}(x' \mid y) + \nabla \log Df_{X \mid Y}(x \mid y)^{\top}(x'-x) \\ \leq& \frac{L'+\eta^{-1}}{2} \|x-x'\|^2.
		\end{aligned}
		\]
		By the assumption on~$\eta$, $a' + \eta^{-1} \geq L' - a'$ and $a' + \eta^{-1} \geq (L' + \eta^{-1})/2$, so the conditional potential is convex, and its condition number is at most $2$. Corollary~2.2 of \cite{qin2026randomizedmala} gives the desired spectral gap bound for $Q_y$.
	\end{proof}
	
	\begin{proposition} \label{pro:location-semi}
		Under the location-mixture assumptions in Section~\ref{sssec:proximal-location}, suppose that $R:=\sup_{z,z'}\|m(z)-m(z')\|<\infty$. Then~\eqref{ine:semi} holds with
		\[
		L'=L,\qquad a'=a-\frac{L^2R^2}{4}.
		\]
	\end{proposition}
	
	The proof of Proposition~\ref{pro:location-semi} is found in Appendix~\ref{app:location-semi}.
	
	Consider the example where a finite number of $z$'s are arranged along a path.
	
	\begin{proposition} \label{pro:hybrid-location}
		Assume the location-mixture conditions in Section~\ref{sssec:proximal-location}. Suppose that $\Z=[p]$, where $p\geq2$ is an integer, $f_Z(z)=1/p$ for $z\in[p]$, and $\|m(z)-m(z+1)\|\leq\kappa$ for $z\in[p-1]$, where $\kappa\geq0$. Choose
		\[
		\eta=\frac{1}{L(1+p^2L\kappa^2)},\qquad
		h_{\max}=\frac{1}{L(2+p^2L\kappa^2)\sqrt{d(1+\log d)}}.
		\]
		Then, for a universal constant $c\in(0,\infty)$, the hybrid kernel satisfies
		\[
		\mathrm{Gap}(H)\geq
		\frac{ca}{Lp^3(1+p^2 L \kappa^2)\sqrt{d(1+\log d)}}e^{-L\kappa^2/[4(1+L\eta)]}.
		\]
	\end{proposition}
	
	\begin{proof}
		The triangle inequality gives $\sup_{z,z'}\|m(z)-m(z')\|\leq(p-1)\kappa\leq p\kappa$.
		Proposition~\ref{pro:location-semi} therefore gives~\eqref{ine:semi} with $L'=L$ and $a'=a-L^2p^2\kappa^2/4$.
		Our choice of $\eta$ satisfies $\eta^{-1}\geq L'-2a'$, so Proposition~\ref{pro:proximal-semi} gives
		$\mathrm{Gap}(Q_y)\geq c/\sqrt{d(1+\log d)}$ uniformly in $y$.
		By the bound in Proposition~\ref{pro:proximal-location}(ii),
		\[
		\begin{aligned}
			\mathrm{Gap}(K)&\geq\frac{a\eta}{1+a\eta}\frac2p
			\left[1-\cos\left(\frac\pi p\right)\right]e^{-L\kappa^2/[4(1+L\eta)]}\\
			&\geq\frac{2a}{Lp^3(1+p^2 L \kappa^2)}e^{-L\kappa^2/[4(1+L\eta)]},
		\end{aligned}
		\]
		where the second inequality uses $a\eta\leq a/L\leq1$ and $1-\cos(\pi/p)\geq2/p^2$.
		Multiply these bounds using Proposition~\ref{pro:hybrid} and adjust the universal constant $c$ to arrive at the desired conclusion.
	\end{proof}
	
	\begin{corollary} \label{cor:mixing-hybrid-location}
		Let $\epsilon > 0$, and let $\mu\in L_*^2(f_X)$.
		Under the assumptions of~\ref{pro:hybrid-location}, as $L/a, p, L \kappa^2, d, \|\mu - f_X\|_{f_X}/\epsilon$ tend to infinity,
		\[
		t(\epsilon;\mu,H)=O\!\left(
		\frac{L}{a}p^3(1+ p^2 L \kappa^2)\sqrt{d(1+\log d)} e^{L\kappa^2/[4(1+L\eta)]}
		\log (\|\mu-f_X\|_{f_X}/\epsilon) \right).
		\]
	\end{corollary}
	\begin{proof}
		The conditional kernels $Q_y$ are positive semi-definite by laziness, and hence so is $H$.
		Apply Corollary~\ref{cor:mixing-sb} with $b=0$ and $s$ equal to the
		reciprocal of the gap lower bound in Proposition~\ref{pro:hybrid-location}.
	\end{proof}
	
	%For comparison, direct Gaussian random-walk Metropolis, with any fixed proposal multiplier $\varsigma>0$, has the bound
	%\[
	%t(\epsilon;\mu,K)=O\!\left(\frac{L}{a}dp^3e^{L\kappa^2/4}
	%\log\|\mu-f_X\|_{f_X}\right)
	%\]
	%from Corollary~\ref{cor:mixing-MH-location}(ii). Both bounds contain the path factor $p^3$, the overlap factor $e^{L\kappa^2/4}$ determined by neighboring means, and $(L/a)\log\|\mu-f_X\|_{f_X}$. The randomized-step MALA implementation replaces the remaining factor $d$ by $(1+Lp^2\kappa^2)\sqrt{d(1+\log d)}$. Thus, for fixed $p$ and $L\kappa^2$, its dimension dependence improves from $d$ to $\sqrt{d\log d}$ as $d\to\infty$. More generally, its bound is smaller in order when
	%\[
	%1+Lp^2\kappa^2=o\!\left(\sqrt{\frac{d}{1+\log d}}\right).
	%\]
	%The extra factor $1+Lp^2\kappa^2$ comes from controlling the curvature of the conditional target and prevents a uniform improvement as the path grows.
	%This compares iteration bounds; MALA additionally requires evaluations of $\nabla\log Df_X$.
	
	\subsection{Random-scan Gibbs sampler} \label{ssec:gibbs}
	
	\subsubsection{Overview}
	
	Partition $\mathbb{R}^d=\prod_{i=1}^B\mathbb{R}^{d_i}$, where $d_i\geq1$ and $\sum_i d_i=d$.
	Let $K$ be the random-scan Gibbs sampler in \ref{ex:gibbs}, with each iteration updating one uniformly chosen block from its full conditional under  $f_X$.
	%Under the mixture representation \eqref{eq:mixture},
	%\[
	% f_{X_i\mid X_{-i}}(\df x_i\mid x_{-i})
	% =\int_{\Z}f_{X_i\mid X_{-i},Z}(\df x_i\mid x_{-i},z)
	%                 f_{Z\mid X_{-i}}(\df z\mid x_{-i}).
	%\]
	For each $z$, take $K_z$ to be the Gibbs sampler with the same blocks and scan probabilities targeting $f_{X\mid Z}(\cdot\mid z)$.
	Proposition~\ref{pro:condition-DA} verifies Condition~\ref{cond:A}.
	
	The component estimate below is an immediate specialization of the approximate tensorization and spectral-gap bound in \citet[Corollary~3.7]{ascolani2024entropy}.
	\begin{proposition} \label{pro:gibbs-component}
		Suppose that $f_{X\mid Z}(\df x\mid z)=\exp[-u_z(x)]\,\df x$, where $u_z$ is differentiable, $a_z$-strongly convex, and has an $L_z$-Lipschitz gradient, for constants $0<a_z\leq L_z<\infty$.
		Then
		\[
		\mathrm{Gap}(K_z)\geq\frac{a_z}{B L_z}.
		\]
		Consequently, if $a_z/L_z\geq a/L$ for $f_Z$-almost every $z$, then
		\[
		\mathrm{Gap}(K) \geq\frac{a}{B L} \mathrm{Gap}(\bar K_0).
		\]
	\end{proposition}
	\begin{proof}
		By \citet[Assumption~A, Lemma~2.4, and Corollary~3.7]{ascolani2024entropy}, these assumptions give a component condition number at most $L_z/a_z$ and the approximate tensorization bound
		\[
		\mathrm{Var}_{f_{X\mid Z}(\cdot\mid z)}(g)
		\leq\frac{L_z}{a_z}\sum_{i=1}^B
		\mathrm{E}\big[\mathrm{Var}(g(X)\mid X_{-i},Z=z)\mid Z=z\big].
		\]
		The sum on the right is $B\mathcal E_{K_z}(g)$.
		The conclusion follows from Corollary \ref{cor:uniform-gap}.
	\end{proof}
	
	More refined block condition numbers from \citet{ascolani2024entropy} may replace $L_z/a_z$ whenever available.
	If these ratios are not uniformly bounded, Proposition~\ref{pro:main} instead yields a weak Poincar\'{e} inequality with $\varphi(z)=B L_z/a_z$.
	
	\subsubsection{Location mixtures}
	
	Assume the location-mixture condition in Section~\ref{sssec:MH-location}.
	Combining Proposition~\ref{pro:gibbs-component} with Proposition~\ref{pro:Kbar0-location} gives the following.
	
	\begin{proposition} \label{pro:gibbs-location}
		Let $\breve K_0$ satisfy the assumptions of Proposition~\ref{pro:Kbar0-location}.
		Then
		\[
		\mathrm{Gap}(K)\geq \frac{a}{BL}\mathrm{Gap}(\breve K_0).
		\]
		In particular:
		\begin{enumerate}
			\item[(i)] If $\sup_{z,z'}\|m(z)-m(z')\|\leq R<\infty$, then
			\[
			\mathrm{Gap}(K)\geq\frac{a}{BL}\exp\left(-\frac{LR^2}{4}\right).
			\]
			\item[(ii)] If $\Z=[p]$, $p\geq2$, $f_Z(z)=1/p$, and $\|m(z)-m(z+1)\|\leq\kappa$ for $z\in[p-1]$, then
			\[
			\mathrm{Gap}(K)\geq\frac{a}{BL}\frac{2}{p}
			\left(1-\cos\frac{\pi}{p}\right)
			\exp\left(-\frac{L\kappa^2}{4}\right).
			\]
		\end{enumerate}
	\end{proposition}
	
	\begin{proof}
		Proposition~\ref{pro:gibbs-component} gives $\mathrm{Gap}(K_z)\geq a/(BL)$. Proposition~\ref{pro:Kbar0-location} gives $\mathrm{Gap}(\bar K_0)\geq\mathrm{Gap}(\breve K_0)$. 
		Combine these bounds in Proposition~\ref{pro:gibbs-component}. The two special cases follow from Corollaries~\ref{cor:Kbar0-location-bounded} and~\ref{cor:Kbar0-location-line}.
	\end{proof}
	
	\begin{corollary} \label{cor:mixing-gibbs-location}
		Let $\epsilon>0$ and let $\mu\in L_*^2(f_X)$.
		Under the respective assumptions of Proposition~\ref{pro:gibbs-location},
		\[
		t(\epsilon;\mu,K)=
		\begin{cases}
			O\bigl((L/a)B e^{LR^2/4} \log(\|\mu-f_X\|_{f_X}/\epsilon)\bigr),
			&\text{in (i)},\\
			O\bigl((L/a)B p^3e^{L\kappa^2/4} \log (\|\mu-f_X\|_{f_X}/\epsilon) \bigr),
			&\text{in (ii)}.
		\end{cases}
		\]
		In (i), the big $O$ notation depicts asymptotics when $\|\mu - f_X\|_{f_X}/\epsilon$, $L/a$, $B$, and $LR^2/4$ tend to infinity;
		in (ii), the big $O$ notation depicts asymptotics when $\|\mu - f_X\|_{f_X}/\epsilon$, $L/a$, $B$, $p$, and $L\kappa^2/4$ tend to infinity.
	\end{corollary}

	\subsubsection{Scale mixtures}
	
	Assume the scale-mixture condition in Section~\ref{sssec:Kbar0-scale}.

	\begin{proposition} \label{pro:gibbs-scale}
		Assume that $\mathrm{Is}(f_{\log Z})>0$.
		Then
		\begin{equation} \nonumber
			\mathrm{Gap}(K)\geq\frac{1}{128B} \frac{a}{L}
			\left[1\wedge\frac{\mathrm{Is}(f_{\log Z})^2}{(1+L/a)d}\right].
		\end{equation}
		As a consequence, for $\epsilon > 0$ and $\mu \in L_*^2(f_X)$, as $B$, $L/a$, $d$, $\|\mu-f_X\|_{f_X}/\epsilon \to \infty$,
		\[
		t(\epsilon; \mu, K) = O \left( \left( \frac{L}{a} \right)^2 B d \log \frac{\|\mu-f_X\|_{f_X}}{\epsilon} \right).
		\]
	\end{proposition}
	\begin{proof}
		Apply Proposition~\ref{pro:gibbs-component} with the scale-invariant ratio $a_z/L_z=a/L$, followed by Proposition~\ref{pro:Kbar0-scale}.
	\end{proof}

	\begin{corollary} \label{cor:mixing-gibbs-scale}
		Let $f_Z$ have the inverse-gamma
		law~\eqref{eq:ig} with $\theta>0$.
		Let $\epsilon>0$ and let $\mu \in L_*^2(f_X)$.
		Then, as $B, d, L/a, \|\mu-f_X\|_{f_X}/\epsilon \to \infty$ and $\theta \downarrow0$,
		\[
		t(\epsilon;\mu,K)
		=O\left(
		\left( \frac{L}{a} \right)^2 B d \theta^{-2} \log\frac{\|\mu-f_X\|_{f_X}}{\epsilon}
		\right).
		\]
	\end{corollary}
	\begin{proof}
		Apply Corollary~\ref{cor:mixing-sb} with $b=0$ and
		\[
		s=\frac{128BL}{a}
		\left[1\vee\frac{(1+L/a)d}{\mathrm{Is}(f_{\log Z})^2}\right],
		\]
		Proposition~\ref{pro:Is-ig} gives $s=O((L/a)^2Bd\theta^{-2})$ as $\theta\downarrow0$.
		The desired result then follows.
	\end{proof}
	
	\subsection{Simulation study} \label{ssec:scale-simulation}
	
	Let $f_Z$ be the inverse-gamma law \eqref{eq:ig} with parameter $\theta > 0$.
	The mixing-time upper bounds in Corollaries~\ref{cor:mixing-MH-scale},
	\ref{cor:mixing-proximal-scale}, and~\ref{cor:mixing-gibbs-scale}
	deteriorate as $\theta\downarrow 0$, with a stronger dependence on
	$\theta$ for \ref{ex:MH} and \ref{ex:proximal} than for \ref{ex:gibbs}.
	The first two bounds account for the probability of scales at which
	the component gap estimates for $K_z$ become small.
	
	We compare the three samplers with $f_{X\mid Z}(\cdot\mid z)=\mathrm{N}_d(0,zI_d)$,
	so $f_X$ is the multivariate
	Student's $t$ distribution with $\theta$ degrees of freedom. 
	Set $d=100$,
	$\sigma^2=\eta=1/(8d)=1/800$, and use $B=d$ one-dimensional blocks for random-scan Gibbs.
	For each $\theta$,
	the three samplers share $n=8{,}192$ independent starting states drawn from
	$\mu=f_X(\cdot\mid \|X\|\le r_{.5})$, where $r_q$ is the target $q$-quantile
	of the radius. This gives $\|\mu-f_X\|_{f_X}=1$ for every $\theta$.
	
	Let $F(r)=f_X(\{x:\|x\|\le r\})$ denote the target radius cumulative
	distribution function (CDF), and let $F_t$ be the radius CDF under
	$\mu K^t$. Set $D_t=\sup_r|F_t(r)-F(r)|$, which is a lower bound on
	$\mathrm{TV}(\mu K^t,f_X)$. We estimate $D_t$ by
	$\widehat D_t=\sup_r|\widehat F_t(r)-F(r)|$, where $\widehat F_t$ is the
	empirical CDF of the radii across chains at iteration $t$.
%	Table~\ref{tab:scale-simulation} reports $\widehat D_t$ at $t=10^5$.
%	If the $4{,}096$ endpoints were independent draws from $f_X$, the
%	Dvoretzky--Kiefer--Wolfowitz inequality \citep{massart1990tight} would give
%	a $95\%$ upper bound of $\sqrt{\log(2/0.05)/(2\cdot 4{,}096)}
%	\approx 0.0212$ for $\widehat D_t$.
%	
%	\begin{table}[H]
%		\centering
%		\begin{tabular}{@{}rrrr@{}}
%			\toprule
%			$\theta$ & \ref{ex:MH} & \ref{ex:proximal} & \ref{ex:gibbs} \\
%			\midrule
%			$0.05$ & $0.5001$ & $0.5001$ & $0.4393$ \\
%			$0.10$ & $0.5001$ & $0.5002$ & $0.3820$ \\
%			$0.25$ & $0.4848$ & $0.4775$ & $0.2198$ \\
%			$0.50$ & $0.4201$ & $0.3792$ & $0.0963$ \\
%			$1.00$ & $0.3056$ & $0.2385$ & $0.0129$ \\
%			$2.00$ & $0.1656$ & $0.1031$ & $0.0117$ \\
%			$4.00$ & $0.0651$ & $0.0242$ & $0.0154$ \\
%			\bottomrule
%		\end{tabular}
%		\caption{Empirical radial Kolmogorov--Smirnov distances in dimension $100$
%		after $10^5$ iterations.}
%		\label{tab:scale-simulation}
%	\end{table}
%	
%	The radial discrepancies in Table~\ref{tab:scale-simulation} generally
%	increase as $\theta$ decreases. At $\theta=0.25$, the discrepancies for
%	\ref{ex:MH} and \ref{ex:proximal} remain close to the initial population
%	discrepancy of $0.5$, while that for \ref{ex:gibbs} is $0.2198$.

%	\paragraph{Radial mixing times.}
	We ran an experiment at $\theta\in\{0.15,0.32,0.65,0.92,1.3,2.6,5.2\}$.
	We used $m=24$ prespecified checkpoints,
	including iteration zero, with successive positive checkpoints spaced
	by a factor of about $\sqrt{2}$ from $128$ through $2^{18}=262{,}144$.
	The radial mixing time on this grid is the first checkpoint with
	$D_t\le\epsilon$, where $\epsilon=0.30$. We estimate it by the first
	checkpoint with $\widehat D_t\le\epsilon$.
	
\begin{figure}[!htbp]
		\centering
		\includegraphics[width=\linewidth]{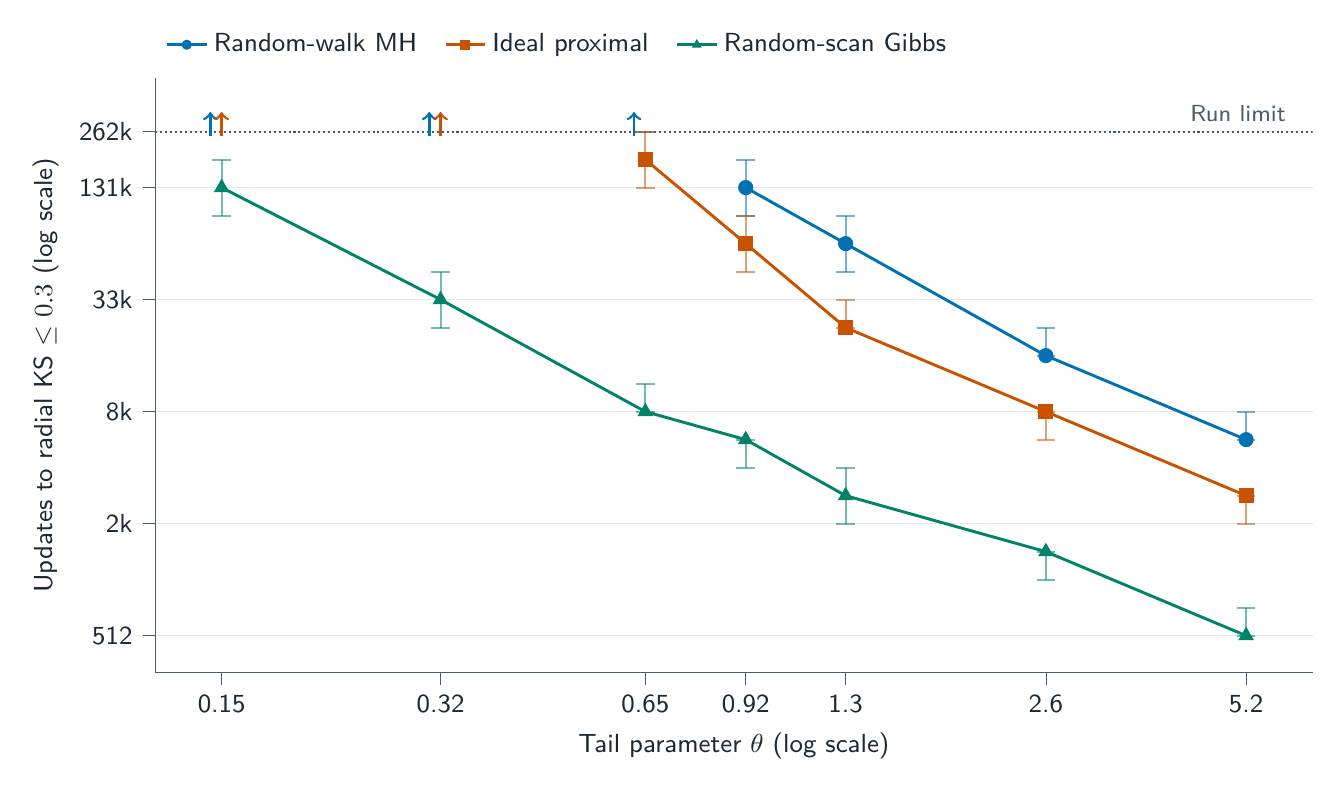}
		\caption{Estimated radial mixing times at tolerance $0.30$.
		Arrows at the run limit indicate settings in which the empirical
		distance exceeded $0.30$ at every checkpoint through $2^{18}$ iterations.}
		\label{fig:radial-mixing-time}
	\end{figure}

	At each checkpoint the radii are independent across chains. The
	Dvoretzky--Kiefer--Wolfowitz inequality \citep{massart1990tight}
	bounds the probability that
	$\sup_r|\widehat F_t(r)-F_t(r)|>b$ by $2\exp(-2nb^2)$. Since this
	supremum bounds $|\widehat D_t-D_t|$, a union bound over the $m$
	checkpoints gives
	\[
	  |\widehat D_t-D_t|\le b_n
	  \quad\text{at every checkpoint},\qquad
	  b_n=\sqrt{\frac{\log(2m/0.05)}{2n}}\approx0.02047
	\]
	with probability at least $95\%$ for each sampler and $\theta$.
	The simultaneous bound also covers stopping at the first checkpoint
	with $\widehat D_t+b_n\le\epsilon$.
	In Figure~\ref{fig:radial-mixing-time}, the upper and lower ends of each error bar are the first checkpoints with
	$\widehat D_t+b_n\le\epsilon$ and with $\widehat D_t-b_n\le\epsilon$, respectively.
	If there is no checkpoint satisfying the upper-end condition, the
	confidence interval has no finite upper endpoint.

	At $\theta=0.15$, both \ref{ex:MH} and \ref{ex:proximal} ran to the
	limit, with $\widehat D_t\ge 0.4972$ at all monitored checkpoints.
	\ref{ex:gibbs} first crossed at $131{,}072$ iterations, with a $95\%$
	confidence interval of $[92{,}682,185{,}364]$ on the checkpoint grid.
	The separation at small $\theta$
	is consistent with the stronger dependence on $\theta$ in the Metropolis
	and proximal bounds in Corollaries~\ref{cor:mixing-MH-scale}
	and~\ref{cor:mixing-proximal-scale}, relative to the Gibbs bound in
	Corollary~\ref{cor:mixing-gibbs-scale}.

	\section*{Acknowledgment}
	The author thanks Eric Chi and Guanyang Wang for the helpful discussions that inspired this project. 
	
	\appendix

	\section{Basic Operator Theory for Reversible Markov Chains} \label{app:operator}
	
	In this subsection, we lay out the basic operator theory for Markov chains that will be used in some of our proofs.
	
	Let $(\Omega, \mathcal{F}, \varpi)$ be a probability space.
	Recall that the space $L_0^2(\varpi)$ consists of $g: \Omega \to \mathbb{R}$ such that $\varpi g = 0$ and that $\|g\|_{\varpi} = \sqrt{\varpi (g^2)} < \infty$.
	Two functions in $L_0^2(\varpi)$ are considered equal if they are equal $\varpi$-almost everywhere.
	For $g_1, g_2 \in L_0^2(\varpi)$ and $c \in \mathbb{R}$, one may define $(g_1 + g_2)(\cdot) = g_1(\cdot) + g_2(\cdot)$ and $(cg_1)(\cdot) = c g_1(\cdot)$, thus making $L_0^2(\varpi)$ a vector space.
	One may further endow the space with an inner product of the form $\langle g_1, g_2 \rangle_{\varpi} = \int_{\Omega} g_1(\omega) g_2(\omega) \, \varpi(\df \omega)$.
	It turns out that this makes $L_0^2(\varpi)$ a Hilbert space.
	For $g \in L_0^2(\varpi)$, the quantity $\|g\|_{\varpi}$ can be seen as the norm $\sqrt{\langle g, g \rangle_{\varpi}}$.
	
	Let $(\Omega_1, \mathcal{F}_1, \varpi_1)$ and $(\Omega_2, \mathcal{F}_2, \varpi_2)$ be two probability spaces.
	Let $T: \Omega_1 \times \mathcal{F}_2 \to [0,1]$ be an Mtk such that $\varpi_1 T(\cdot) := \int_{\Omega_1} \varpi_1(\df \omega) \, T(\omega, \cdot)  = \varpi_2(\cdot)$.
	It defines a linear transformation from $L_0^2(\varpi_2)$ to $L_0^2(\varpi_1)$ in the following way: For $g \in L_0^2(\varpi_2)$, $T g(\cdot) = \int_{\Omega_2} T(\cdot, \df \omega) \, g(\omega)$.
	By Jensen's inequality, the operator is bounded in the sense that
	\[
	\|T\|_{\varpi_2 \to \varpi_1} := \sup_{g \in L_0^2(\varpi_2) \setminus \{0\}} \frac{\|Tg\|_{\varpi_1}}{\|g\|_{\varpi_2}} \leq 1.
	\]
	If $\varpi_1 = \varpi_2$, we simply write $\|T\|_{\varpi_2 \to \varpi_1}$ as $\|T\|_{\varpi_1}$.
	
	Let $A: L_0^2(\varpi_2) \to L_0^2(\varpi_1)$ be a bounded linear transformation.
	The adjoint of~$A$ is a bounded linear transformation $A^{\dagger}: L_0^2(\varpi_1) \to L_0^2(\varpi_2)$ such that, for $g_1 \in L_0^2(\varpi_1)$ and $g_2 \in L_0^2(\varpi_2)$, $\langle g_1, A g_2 \rangle_{\varpi_1} = \langle A^{\dagger} g_1, g_2 \rangle_{\varpi_2}$.
	It holds that $(A^{\dagger})^{\dagger} = A$.
	If $A = A^{\dagger}$, we say $A$ is self-adjoint.
	When $A: L_0^2(\varpi_1) \to L_0^2(\varpi_1)$ corresponds to an Mtk that is reversible with respect to $\varpi_1$, it is self-adjoint if and only if the Mtk is reversible with respect to~$\varpi_1$.
	If $A: L_0^2(\varpi_1) \to L_0^2(\varpi_1)$ is self-adjoint, and $\langle g, A g\rangle_{\varpi_1} \geq 0$ for each $g \in L_0^2(\varpi_1)$, we say $A$ is positive semi-definite.
	
	The following is a useful lemma that can be easily established via the spectral decomposition of self-adjoint operators.
	\begin{lemma} \label{lem:spectral}
		Let $K(\cdot, \cdot)$ be an Mtk that is reversible with respect to a probability distribution~$\varpi$, so that $K$ is a linear operator on $L_0^2(\varpi)$.
		Assume that a polynomial function $u \mapsto h(u) = \sum_{j=0}^m c_j u^j$ is non-negative on $[-1,1]$.
		Then, for $g \in L_0^2(\varpi)$, $\langle g, h(K) g \rangle_{\varpi} \geq 0$, where $h(K) = \sum_{j=0}^m c_j K^j$.
		When $K$ is positive semi-definite, the result continues to hold if $h(\cdot)$ is non-negative on $[0,1]$ rather than $[-1,1]$.
	\end{lemma}

	\section{Basic properties of DA chains} \label{app:da}

	For two random elements $(U,V)$ taking values in the Polish space $\mathsf{U} \times \mathsf{V}$, we can always define a DA Mtk
	\[
	K_{U,V}(u, \df u') = \int_{\mathsf{V}} f_{V \mid U}(\df v \mid u) \, f_{U \mid V}(\df u' \mid v),
	\]
	which is reversible with respect to~$f_U$, the marginal distribution of~$U$.
	We briefly mention a few existing results regarding the convergence properties of $K_{U,V}$.
	
	\begin{proposition} 
		The linear operator $K_{U,V}: L_0^2(f_U) \to L_0^2(f_U)$ is positive semi-definite.
	\end{proposition}
	
	\begin{proposition} \label{pro:max-corr} \citep{liu1994covariance}
		It holds that
		\[
		\mathrm{Gap}(K_{U,V}) = 1-\left[\sup_{g,h} \frac{\mathrm{Cov} [g(U), h(V)] }{\sqrt{\mathrm{Var}[g(U)] \, \mathrm{Var}[h(V)] }}\right]^2,
		\]
		where the supremum is over measurable functions with finite, strictly positive variances. Set the maximal correlation to zero if either centered $L^2$ space is trivial.
	\end{proposition}
	\begin{proof}
		Let $Ah(u)=\mathbb E[h(V)\mid U=u]$ on the centered $L^2$ spaces. Its adjoint is $A^\dagger g(v)=\mathbb E[g(U)\mid V=v]$, and $K_{U,V}=AA^\dagger$. The displayed supremum equals $\|A\|$, whereas $\|AA^\dagger\|=\|A\|^2$. Positivity of $K_{U,V}$ then gives the identity.
	\end{proof}

	Proposition~\ref{pro:max-corr} indicates that the convergence property of $K_{U,V}$ is largely determined by how much the distribution of~$V$ depends on~$U$.
	Smaller dependence implies faster convergence.

	Finally, Lemma~\ref{lem:da-elements} below allows one to simplify Dirichlet forms associated with $K_{U,V}$ in certain settings.
	
	\begin{lemma} \label{lem:da-elements}
		Suppose that there is a $\sigma$-finite measure $\nu$ such that $f_{V \mid U}(\cdot \mid u) \ll \nu(\cdot)$ for $f_U$-almost every $u \in \mathsf{U}$.
		For $f_U$-almost every $u \in \mathsf{U}$ and every measurable $\mathsf{A} \subset \mathsf{U}$,
		\[
		K_{U,V}(u, \mathsf{A}) \geq \int_{\mathsf{A}} \left[ \int_{\mathsf{V}} \sqrt{\frac{\df f_{V \mid U}(\cdot \mid u)}{\df \nu(\cdot)} \, \frac{\df f_{V \mid U}(\cdot \mid u')}{\df \nu(\cdot)} } \, \df \nu \right]^2  f_U(\df u').
		\]
	\end{lemma}
	
	\begin{proof}
		For $f_U$-almost every $u \in \mathsf{U}$ and measurable $\mathsf{A} \subset \mathsf{U}$,
		\[
		K_{U,V}(u, \mathsf{A}) = \int_{\mathsf{A}} \int_{\mathsf{V}} \frac{\df f_{V \mid U}(\cdot \mid u)}{\df \nu(\cdot)}  \, \frac{\df f_{V \mid U}(\cdot \mid u')/\df \nu(\cdot)}{\df f_V / \df \nu} \, \df \nu \, f_U(\df u').
		\]
		By Jensen's inequality,
		\[
		\begin{aligned}
			\int_{\mathsf{V}} \frac{\df f_{V \mid U}(\cdot \mid u)}{\df \nu(\cdot)}  \, \frac{\df f_{V \mid U}(\cdot \mid u')/\df \nu(\cdot)}{\df f_V / \df \nu} \, \df \nu &= \int_{\mathsf{V}} \frac{\df f_{V \mid U}(\cdot \mid u)/\df \nu(\cdot)}{\df f_V/\df \nu}  \, \frac{\df f_{V \mid U}(\cdot \mid u')/\df \nu(\cdot)}{\df f_V / \df \nu} \, \frac{\df f_V}{\df \nu} \, \df \nu\\
			&\geq \left[ \int_{\mathsf{V}} \sqrt{\frac{\df f_{V \mid U}(\cdot \mid u)/\df \nu(\cdot) \; \df f_{V \mid U}(\cdot \mid u')/\df \nu(\cdot)}{(\df f_V /\df \nu)^2} }  \, \frac{\df f_V}{\df \nu} \, \df \nu \right]^2 \\
			&= \left[ \int_{\mathsf{V}} \sqrt{\frac{\df f_{V \mid U}(\cdot \mid u)}{\df \nu(\cdot)}  \, \frac{\df f_{V \mid U}(\cdot \mid u')}{\df \nu(\cdot)}} \, \df \nu \right]^2.
		\end{aligned}
		\]
		The desired result then follows.
	\end{proof}

	\section{Proofs} \label{app:proof}
	
	\subsection{Proof of Proposition \ref{pro:convergence-Poincare}} \label{app:convergence-Poincare}
	
	Before stating the proof, we note that, by Lemma~\ref{lem:spectral} in Appendix~\ref{app:operator}, for the positive semi-definite Mtk~$T$ and any $g \in L_0^2(\varpi)$, $\langle g, T g \rangle_{\varpi} \geq \langle g, T^2 g \rangle_{\varpi}$.
	Thus, $\mathcal{E}_T(g) \leq \mathcal{E}_{T^2}(g)$.
	In particular, under the assumptions of Proposition~\ref{pro:convergence-Poincare},~$T^2$ satisfies an $(s,b)$-Poincar\'{e} inequality.

	\begin{proof}
		Let $\mu$ be a distribution in $L_*^2(\varpi)$, and let $t$ be a positive integer.
		
		By Lemma 22.1.3 of \cite{douc2018markov}, $\mu T^t \in L_*^2(\varpi)$.
		By the Cauchy-Schwarz inequality and the reversibility of~$T$,
		\begin{equation}
			\begin{aligned}\label{ine:muTt}
				\mathrm{TV}(\mu T^t , \varpi) 
				&= \sup_{g: \,\Omega \to [0,1]} \left| \mu T^t g - \varpi g \right| \\
				&= \sup_{\varpi g = 0, \; \|g\|_{\mathrm{osc}} = 1} \left| \mu T^t g - \varpi g \right| \\
				&= \sup_{g \in L_0^2(\varpi), \; \|g\|_{\mathrm{osc}} = 1} \left| \mu T^t g - \varpi T^t g \right| \\
				&= \sup_{g \in L_0^2(\varpi), \; \|g\|_{\mathrm{osc}} = 1} \left| \left\langle \frac{\df \mu}{\df \varpi} - 1, T^t g \right\rangle_{\varpi} \right| \\
				& \leq \|\mu - \varpi\|_{\varpi} \, \sup_{g \in L_0^2(\varpi), \; \|g\|_{\mathrm{osc}} = 1} \|T^t g\|_{\varpi}.
			\end{aligned}
		\end{equation}
		By the $(s,b)$-Poincar\'{e} inequality for $T^2$ and the self-adjoint-ness of~$T$, for any $g \in L_0^2(\varpi)$ and positive integer~$t'$,
		\[
		\begin{aligned}
			\|T^{t'-1} g\|_{\varpi}^2 &\leq s \left( \|T^{t'-1} g\|_{\varpi}^2 - \langle T^{t'-1} g, T^2 T^{t'-1} g \rangle_{\varpi} \right) + b\|T^{t'-1} g\|_{\mathrm{osc}}^2 \\
			&\leq s \left( \|T^{t'-1} g\|_{\varpi}^2 - \|T^{t'} g \|_{\varpi}^2 \right) + b\| g\|_{\mathrm{osc}}^2.
		\end{aligned}
		\]
		(The function $T^0 g$ is $g$ by convention.)
		Then
		\[
		\|T^{t'} g\|_{\varpi}^2 \leq \left( 1 - \frac{1}{s} \right) \|T^{t'-1} g\|_{\varpi}^2 + \frac{b}{s} \|g\|_{\mathrm{osc}}^2.
		\]
		By induction, for $g \in L_0^2(\varpi)$ and positive integer~$t$,
		\[
		\begin{aligned}
			\|T^t g\|_{\varpi}^2 &\leq \left(1-\frac{1}{s}\right)^t \|g\|_{\varpi}^2 + b \left[1-\left(1- \frac{1}{s} \right)^t \right] \, \left\| g \right\|_{\mathrm{osc}}^2 \\
			&\leq \left[ \frac{1}{4} \left(1-\frac{1}{s}\right)^t + b\left\{1-\left(1-\frac{1}{s}\right)^t\right\} \right] \left\| g \right\|_{\mathrm{osc}}^2,
		\end{aligned}
		\]
		where the second line follows from Popoviciu's inequality.
		The desired bound then follows from~\eqref{ine:muTt}.
	\end{proof}

	\subsection{Proof of Proposition~\ref{pro:condition-DA}} \label{app:condition-DA}
	
	\begin{proof}
		Let $g: \X \to \mathbb{R}$ be such that $\|g\|_{f_X}^2 < \infty$.
		It is easy to see that
		\[
		\begin{aligned}
			\mathcal{E}_K(g) &= \|g\|_{f_X}^2 - \int_{\mathsf{Y}} f_Y(\df y) \, \left[ \int_{\X} g(x) \, f_{X \mid Y}(\df x \mid y) \right]^2.
		\end{aligned}
		\]
		Let $(X,Y,Z)$ have the joint distribution $f_Z(\df z) \, f_{X \mid Z}(\df x \mid z) \, f_{Y \mid X}(\df y \mid x)$.
		Then
		\[
		\begin{aligned}
			\int_{\Z} f_Z(\df z) \, \mathcal{E}_{K_z}(g) &= \frac{1}{2} \int f_Z(\df z) \, f_{X \mid Z}(\df x \mid z) \, f_{Y \mid X}(\df y \mid x) \, f_{X \mid Y, Z}(\df x' \mid y, z) \, [g(x) - g(x')]^2 \\
			&= \frac{1}{2} \int f_{X,Y,Z}(\df x, \df y, \df z) \, f_{X \mid Y, Z}(\df x' \mid y, z) \, [g(x) - g(x')]^2 \\
			&= \|g\|_{f_X}^2 - \int_{\Y \times \Z} f_{Y,Z}(\df y, \df z) \, \left[ \int_{\X} g(x) \, f_{X \mid Y,Z}(\df x \mid y,z) \right]^2 \\
			&= \|g\|_{f_X}^2 - \int_{\Y} f_Y(\df y) \int_{\Z} f_{Z \mid Y}(\df z \mid y) \, \left[ \int_{\X} g(x) \, f_{X \mid Y,Z}(\df x \mid y,z) \right]^2.
		\end{aligned}
		\]
		It then follows from Jensen's inequality that $\mathcal{E}_K(g) \geq \int_{\Z} f_Z(\df z) \, \mathcal{E}_{K_z}(g)$.
	\end{proof}

	\subsection{Proof of Theorem \ref{thm:main}} \label{app:main}
	
	\subsubsection{Overview}
	
	To prove Theorem \ref{thm:main}, we introduce a few new Mtks that will facilitate our analysis.
	First, define
	\[
	\bar{K}_2(z, \df z') = \int_{\X^2} f_{X \mid Z}(\df x \mid z) \, K^2(x, \df x') \, f_{Z \mid X}(\df z' \mid x').
	\]
	When $t \in \{0,2\}$, define the following analogue of $\bar{K}_t$ on $(\X,\mathcal{B})$:
	\[
	\bar{K}_t^*(x, \df x''') = \int_{\X\times\Z\times\X} K^{t/2}(x, \df x') \, f_{Z \mid X}(\df z \mid x') \, f_{X \mid Z}(\df x'' \mid z) \, K^{t/2}(x'', \df x''').
	\]
	The Mtk $\bar{K}_t^*$ is reversible with respect to $f_X$.
	
	The following proposition is useful for connecting the $\bar{K}_t$'s and $\bar{K}^*_t$'s.
	
	\begin{proposition} \label{pro:conjugate-1} 
		Each of the following statements holds.
		\begin{enumerate}
			\item[(i)] If $t \in \{0,2\}$, then $\mathrm{Gap}(\bar{K}_t) = \mathrm{Gap}(\bar{K}^*_t)$.
			\item[(ii)] If $t \in \{0,2\}$, and $\bar{K}_t$ satisfies an $(s,b)$-Poincar\'{e} inequality for some $s, b \geq 0$, then $\bar{K}_t^*$ satisfies an $(s+1, b)$-Poincar\'{e} inequality.
			\item[(iii)] If $K$ is positive semi-definite and $t,t'\in\{0,1,2\}$ with $t<t'$, then any $(s,b)$-Poincar\'{e} inequality satisfied by $\bar{K}_t$ is also satisfied by $\bar{K}_{t'}$ for the same values of $(s,b)$.
		\end{enumerate}
	\end{proposition}
	
	Proposition~\ref{pro:conjugate-1} is largely similar to a technical result in \citep{gao2026weak}.
	We provide a proof in Appendix~\ref{app:conjugate} for completeness.
	
	Theorem~\ref{thm:main}(i) follows from Proposition~\ref{pro:conjugate-1}(ii) and Lemma~\ref{lem:main}(i). In part (ii), positivity of $K$ permits the transfer from $\bar K_1$ to $\bar K_2$ in Proposition~\ref{pro:conjugate-1}(iii); applying (ii) of that proposition and Lemma~\ref{lem:main}(ii) proves the claim.
	
	\begin{lemma} \label{lem:main}
		Assume that Condition \ref{cond:A} holds.
		Suppose that there exists $s \geq 0$ such that, for $f_Z$-almost every $z \in \Z$, the Mtk $K_z$ satisfies an $(s,b_z)$-Poincar\'{e} inequality, where $b_z \geq 0$, and $z \mapsto b_z$ is measurable.
		Then each of the following statements holds:
		\begin{enumerate}
			\item [(i)] If $\bar{K}_0^*$ satisfies an $(s',b')$-Poincar\'{e} inequality for some non-negative numbers $s'$ and $b'$, then $K$ satisfies an $(s'',b'')$-Poincar\'{e} inequality, where $s'' = s' s$, and $b'' = b' + s' \int_{\Z} b_z \, f_Z(\df z) $.
			\item [(ii)] If $\bar{K}_2^*$ satisfies an $(s',b')$-Poincar\'{e} inequality for some non-negative numbers $s'$ and $b'$, then $K$ satisfies an $(s'',b'')$-Poincar\'{e} inequality, where $s'' = s'(s+2)$, and $b'' = b' + s' \int_{\Z} b_z \, f_Z(\df z) $.
		\end{enumerate}
	\end{lemma}
	
	It remains to prove Lemma~\ref{lem:main}.
	
	\subsubsection{Proof for Part (i) of Lemma~\ref{lem:main}}
	\begin{proof}
		Assume that $\bar{K}_0^*$ satisfies an $(s',b')$-Poincar\'{e} inequality, where $s', b' \geq 0$.
		Let $g \in L_0^2(f_X)$ be arbitrary.
		Then
		\begin{equation} \label{ine:poincare-kbarstar0}
			\|g\|_{f_X}^2 \leq s' \, \mathcal{E}_{\bar{K}_0^*}(g) + b' \|g\|_{\mathrm{osc}}^2 .
		\end{equation}
		Now, 
		\begin{equation} \label{eq:Dirichlet-Kbarstar0}
			\begin{aligned}
				\mathcal{E}_{\bar{K}_0^*}(g) &= \|g\|_{f_X}^2 - \int f_X(\df x) \, g(x) \, f_{Z \mid X}(\df z \mid x) \, f_{X \mid Z}(\df x' \mid z) \, g(x') \\
				&= \int_{\Z} f_Z(\df z) \int_{\X} f_{X \mid Z}(\df x \mid z) \, g(x)^2 - \int_{\Z} f_Z(\df z)  \left[ \int_{\X} g(x) \, f_{X \mid Z}(\df x \mid z) \right]^2 \\
				&= \int_{\Z} f_Z(\df z) \, \|g_z\|_{f_{X \mid Z}(\cdot \mid z)}^2,
			\end{aligned}
		\end{equation}
		where, for $f_Z$-almost every $z \in \Z$,
		\[
		g_z(\cdot) = g(\cdot) - \int_{\X} g(x) \, f_{X \mid Z}(\df x \mid z) \in L_0^2(f_{X \mid Z}(\cdot \mid z)).
		\]
		To proceed, apply the Poincar\'{e} inequality for $K_z$ to the function $g_z$, and obtain
		\[
		\begin{aligned}
			\| g_z \|_{f_{X \mid Z}(\cdot \mid z)}^2   \leq& s \, \mathcal{E}_{K_z}(g_z) + b_z \|g_z\|_{\mathrm{osc}}^2 \\
			=& \frac{s}{2} \int_{\X^2} f_{X \mid Z}(\df x \mid z) \, K_z(x, \df x') \, [g(x) - g(x')]^2 + b_z\|g_z\|_{\mathrm{osc}}^2  \\
			\leq & s \, \mathcal{E}_{K_z}(g) + b_z \|g\|_{\mathrm{osc}}^2,
		\end{aligned}
		\]
		where the last inequality holds for $f_Z$-almost every~$z$.
		Then, combining \eqref{eq:Dirichlet-Kbarstar0} with the most recent display gives
		\begin{equation*}
			\begin{aligned}
				\mathcal{E}_{\bar{K}_0^*}(g) &\leq s  \int_{\Z} f_Z(\df z) \, \mathcal{E}_{K_z}(g) + \|g\|_{\mathrm{osc}}^2 \int_{\Z} f_Z(\df z) \, b_z .
			\end{aligned}
		\end{equation*}
		It then follows from \eqref{ine:poincare-kbarstar0} and Condition~\ref{cond:A} that
		\begin{equation} \nonumber
			\begin{aligned}
				\|g\|_{f_X}^2
				\leq& s's \int_{\Z} f_Z(\df z) \, \mathcal{E}_{K_z}(g) + \left[ b' + s' \int_{\Z} f_Z(\df z) \, b_z \right] \|g\|_{\mathrm{osc}}^2 \\
				\leq &  s' s \, \mathcal{E}_K(g) + \left[ b' + s' \int_{\Z} f_Z(\df z) \, b_z \right] \|g\|_{\mathrm{osc}}^2.
			\end{aligned}
		\end{equation}
		This completes the proof.
	\end{proof}
	
	\subsubsection{Proof for Part (ii) of Lemma~\ref{lem:main}}
	\begin{proof}
		Assume that $\bar{K}_2^*$ satisfies an $(s',b')$-Poincar\'{e} inequality.
		Let $g \in L_0^2(f_X)$ be arbitrary.
		Then 
		\begin{equation} \label{ine:poincare-kbarstar}
			\|g\|_{f_X}^2 \leq s' \, \mathcal{E}_{\bar{K}_2^*}(g) + b' \|g\|_{\mathrm{osc}}^2 .
		\end{equation}
		Now,
		\begin{equation} \label{eq:Dirichlet-Kbarstar}
			\begin{aligned}
				\mathcal E_{\bar K_2^*}(g)
				&=\|g\|_{f_X}^2-\int_{\Z}f_Z(\df z)
				\left[\int_{\X}Kg(x)f_{X\mid Z}(\df x\mid z)\right]^2\\
				&=\|g\|_{f_X}^2-\|Kg\|_{f_X}^2
				+\int_{\Z}\|(Kg)_z\|_{f_{X\mid Z}(\cdot\mid z)}^2 f_Z(\df z),
			\end{aligned}
		\end{equation}
		where, for $f_Z$-almost every $z \in \Z$,
		\[
		(Kg)_z(\cdot) = Kg(\cdot) - \int_{\X} Kg(x) \, f_{X \mid Z}(\df x \mid z) \in L_0^2(f_{X \mid Z}(\cdot \mid z)).
		\]
		To proceed, apply the Poincar\'{e} inequality for $K_z$ to the function $(Kg)_z$, and obtain
		\[
		\begin{aligned}
			&\| (Kg)_z \|_{f_{X \mid Z}(\cdot \mid z)}^2  \\ \leq& \, s \, \mathcal{E}_{K_z}((Kg)_z) + b_z \|(Kg)_z\|_{\mathrm{osc}}^2 \\
			=& \, \frac{s}{2}  \int_{\X^2} f_{X \mid Z}(\df x \mid z) \, K_z(x, \df x') \, [Kg(x) - Kg(x')]^2 + b_z \|(Kg)_z\|_{\mathrm{osc}}^2 \\
			\leq& \, s \, \mathcal{E}_{K_z}(K g) + b_z\|Kg\|_{\mathrm{osc}}^2,
		\end{aligned}
		\]
		where the final inequality holds for $f_Z$-almost every~$z$.
		Furthermore, $\|Kg\|_{\mathrm{osc}}$ is upper bounded by $\|g\|_{\mathrm{osc}}$.
		Then, combining \eqref{eq:Dirichlet-Kbarstar} with the most recent display, and applying Condition~\ref{cond:A} yields
		\begin{equation*}
			\begin{aligned}
				\mathcal{E}_{\bar{K}_2^*}(g) &\leq \|g\|_{f_X}^2 - \|Kg\|_{f_X}^2 + s \int_{\Z} f_Z(\df z) \,  \mathcal{E}_{K_z}(Kg) + \|g\|_{\mathrm{osc}}^2 \int_{\Z} f_Z(\df z) \, b_z . \\
				&\leq \|g\|_{f_X}^2 + ( s-1) \|Kg\|_{f_X}^2 - s \langle Kg, K^2 g \rangle_{f_X} + \|g\|_{\mathrm{osc}}^2 \int_{\Z} f_Z(\df z) \, b_z.
			\end{aligned}
		\end{equation*}
		It then follows from \eqref{ine:poincare-kbarstar} that
		\begin{equation} \label{ine:K3}
			\begin{aligned}
				& s' s \langle g, K^3 g \rangle_{f_X} - s'( s-1) \langle g, K^2 g \rangle_{f_X}- (s'-1) \|g\|_{f_X}^2 \leq \left[ b' + s' \int_{\Z} f_Z(\df z) \, b_z \right] \|g\|_{\mathrm{osc}}^2.
			\end{aligned}
		\end{equation}
		Here reversibility gives
		\[
		\|Kg\|_{f_X}^2=\langle g,K^2g\rangle_{f_X},\qquad
		\langle Kg,K^2g\rangle_{f_X}=\langle g,K^3g\rangle_{f_X}.
		\]
		
		To proceed, we find a lower bound on the left-hand-side of \eqref{ine:K3} in terms of $\langle g, K g \rangle_{f_X}$ and $\|g\|_{f_X}^2$.
		For $u \in [-1,1]$, it holds that
		\[
		s's u^3 - s'(s-1)u^2 - s'(s+2) u + s'(s+1) = s' (u-1)(su^2 + u - s-1) \geq 0,
		\]
		which is equivalent to
		\[
		s' s u^3 - s' (s-1) u^2 - (s'-1) \geq 1 - s'(s+2) + s'(s+2)u.
		\]
		By Lemma~\ref{lem:spectral} in Appendix~\ref{app:operator}, the above display implies that 
		\[
		\begin{aligned}
			& s' s \langle g, K^3 g \rangle_{f_X} - s' (s-1) \langle g, K^2 g\rangle_{f_X}- (s'-1) \|g \|_{f_X}^2 \\
			&\geq [1-s'(s+2)] \|g\|_{f_X}^2 + s'(s+2) \langle g, Kg \rangle_{f_X} \\
			&= \|g\|_{f_X}^2 - s'(s+2) \, \mathcal{E}_K(g).
		\end{aligned}
		\]
		Thus, by \eqref{ine:K3},
		\[
		\|g\|_{f_X}^2 \leq s'(s+2) \, \mathcal{E}_K(g) + \left[ b' + s' \int_{\Z} f_Z(\df z) \, b_z \right] \|g\|_{\mathrm{osc}}^2.
		\]
		This completes the proof.
		
	\end{proof}

	\subsubsection{Proof of Proposition~\ref{pro:conjugate-1}} \label{app:conjugate}
	
	\begin{proof}
		For $t \in \{0,2\}$, let $A_t: L_0^2(f_Z) \to L_0^2(f_X)$ such that, for $h \in L_0^2(f_Z)$,
		\[
		A_t \, h(x) = \int_{\X \times \Z} K^{t/2}(x, \df x') \, f_{Z \mid X}(\df z \mid x') \, h(z), \quad x \in \X.
		\]
		Then the adjoint of $A_t$ is the operator $A_t^{\dagger}: L_0^2(f_X) \to L_0^2(f_Z)$ such that, for $g \in L_0^2(f_X)$,
		\[
		A_t^{\dagger} \, g(z) = \int_{\X^2} f_{X \mid Z}(\df x \mid z) \, K^{t/2}(x, \df x') \, g(x'), \quad z \in \Z.
		\]
		Here, we have used the fact that $K^{t/2}$ is reversible with respect to $f_X$.
		Note that $\bar{K}_t = A_t^{\dagger} A_t$, and $\bar{K}_t^* = A_t A_t^{\dagger}$.
		
		\begin{enumerate}
			\item[(i)] 
			Let $t \in \{0,2\}$.
			Suppose that $\mathrm{Gap}(\bar{K}_t) \geq 1/s$ for some $s > 0$, i.e., for $h \in L_0^2(f_Z)$,
			\[
			\|h\|_{f_Z}^2 \leq s [\|h\|_{f_Z}^2 - \langle h, A_t^{\dagger} A_t h \rangle_{f_Z}].
			\]
			We may take $s\geq1$, since $A_t^{\dagger} A_t$ is positive semi-definite and we use the gap-one convention on trivial centered spaces.
			By Lemma~\ref{lem:poincare-conjugate} stated after this proof, if $s > 1$, then for $g \in L_0^2(f_X)$,
			\[
			\|g\|_{f_X}^2 \leq s [\|g\|_{f_X}^2 - \langle g, A_t A_t^{\dagger} g \rangle_{f_X}].
			\]
			That is, $\mathrm{Gap}(A_t A_t^{\dagger}) \geq 1/s$ if $s > 1$.
			Letting $s \to 1$ shows the inequality in fact holds for $s \geq 1$.
			Thus, $\mathrm{Gap}(A_t A_t^{\dagger}) \geq \mathrm{Gap}(A_t^{\dagger} A_t)$.
			To reverse the inequality, use Lemma~\ref{lem:poincare-conjugate} again and note that $(A_t^{\dagger})^{\dagger} = A_t$.
			\item [(ii)] This also follows from Lemma~\ref{lem:poincare-conjugate}.
			\item [(iii)]
			Suppose that $t,t'\in\{0,1,2\}$ and $t<t'$, and that $\bar{K}_t$ satisfies an $(s,b)$-Poincar\'{e} inequality, where $s,b \geq 0$.
			The Poincar\'{e} inequality states that, for $h \in L_0^2(f_Z)$,
			\[
			\begin{aligned}
				\|h\|_{f_Z}^2 &\leq s \|h\|_{f_Z}^2 - s \int f_Z(\df z) \, f_{X \mid Z}(\df x \mid z) \, K^t(x, \df x') \, f_{Z \mid X}(\df z' \mid x') \, h(z) \, h(z') + b \|h\|_{\mathrm{osc}}^2 \\
				&= s \|h\|_{f_Z}^2 - s \langle A_0 h, K^tA_0 h \rangle_{f_X} + b \|h\|_{\mathrm{osc}}^2.
			\end{aligned}
			\]
			Since $K$ is positive semi-definite, by Lemma~\ref{lem:spectral}, $\langle A_0 h, K^tA_0 h \rangle_{f_X} \geq \langle A_0 h, K^{t'}A_0 h \rangle_{f_X}$.
			Then $\bar{K}_{t'}$ also satisfies the $(s,b)$-Poincar\'{e} inequality.
		\end{enumerate}
	\end{proof}
	
	The following lemma is adopted from \cite{gao2026weak} with some modifications.
	
	\begin{lemma} \label{lem:poincare-conjugate}
		Let $(\Omega_1, \mathcal{F}_1, \varpi_1)$ and $(\Omega_2, \mathcal{F}_2, \varpi_2)$ be probability spaces.
		Let $A: L_0^2(\varpi_2) \to L_0^2(\varpi_1)$ be a bounded linear transformation.
		Assume that, for $g \in L_0^2(\varpi_1)$, it holds that $\|A^{\dagger} g\|_{\mathrm{osc}} \leq \|g\|_{\mathrm{osc}}$.
		Suppose further that there exist $s > 0$ and $b \geq 0$ such that, for $h \in L_0^2(\varpi_2)$, 
		\begin{equation} \label{ine:AdaggerA}
			\|h\|_{\varpi_2}^2 \leq s ( \|h\|_{\varpi_2}^2 - \langle h, A^{\dagger} A h \rangle_{\varpi_2}) + b \|h\|_{\mathrm{osc}}^2.
		\end{equation}
		Then, for $s_* \in ( (\sqrt{s}-1)^2, (\sqrt{s}+1)^2)$ and $g \in L_0^2(\varpi_1)$,
		\begin{equation} \label{ine:AAdagger-0}
			\|g\|_{\varpi_1}^2 \leq \frac{4ss_*}{4ss_* - (s_*+s-1)^2}  ( \|g\|_{\varpi_1}^2 - \langle g, A A^{\dagger} g \rangle_{\varpi_1}) + \frac{4sb}{4ss_* - (s_*+s-1)^2} \|g\|_{\mathrm{osc}}^2.
		\end{equation}
		In particular, when $s_* = s+1$, one has
		\[
		\|g\|_{\varpi_1}^2 \leq (s+1)  ( \|g\|_{\varpi_1}^2 - \langle g, A A^{\dagger} g \rangle_{\varpi_1}) + b \|g\|_{\mathrm{osc}}^2;
		\]
		The preceding bound also holds when $s=0$. When $s > 1$ and $s_* = s-1$, one has
		\[
		\|g\|_{\varpi_1}^2 \leq s  ( \|g\|_{\varpi_1}^2 - \langle g, A A^{\dagger} g \rangle_{\varpi_1}) + \frac{sb}{s-1}\|g\|_{\mathrm{osc}}^2.
		\]
	\end{lemma}

	\begin{proof}
		Let $g \in L_0^2(\varpi_1)$ be arbitrary. If $s=0$, applying~\eqref{ine:AdaggerA} to $A^\dagger g$ yields $\langle g,AA^\dagger g\rangle_{\varpi_1}\leq b\|g\|_{\mathrm{osc}}^2$, proving the claimed $(1,b)$ bound. Assume henceforth $s>0$.
		Applying \eqref{ine:AdaggerA} to $h = A^{\dagger} g$ gives
		\begin{equation} \label{ine:AAdagger-1}
			(1-s) \langle g, A A^{\dagger} g \rangle_{\varpi_1} + s \|A A^{\dagger} g\|_{\varpi_1}^2 \leq b \|A^{\dagger} g \|_{\mathrm{osc}}^2 \leq b \|g\|_{\mathrm{osc}}^2.
		\end{equation}
		For $s_* \in \mathbb{R}$, it holds that
		\[
		\left\| \sqrt{s} A A^{\dagger} g - \frac{s_*+s-1}{2\sqrt{s}} g \right\|_{\varpi_1}^2 \geq 0.
		\]
		Rearranging, one obtains
		\[
		-\frac{(s_*+s-1)^2}{4s} \|g\|_{\varpi_1}^2 + s_*\langle g, AA^{\dagger} g \rangle_{\varpi_1} \leq (1-s) \langle g, A A^{\dagger} g \rangle_{\varpi_1} + s \|A A^{\dagger} g\|_{\varpi_1}^2.
		\]
		Combined with \eqref{ine:AAdagger-1}, this gives
		\[
		-(s_*+s-1)^2 \|g\|_{\varpi_1}^2 + 4ss_* \langle g, A A^{\dagger} g\rangle_{\varpi_1} \leq 4sb \|g\|_{\mathrm{osc}}^2.
		\]
		Suppose that $s_* \in ( (\sqrt{s}-1)^2, (\sqrt{s}+1)^2)$.
		Then $4ss_* - (s_*+s-1)^2 > 0$, and the most recent display yields~\eqref{ine:AAdagger-0}.
	\end{proof}
	
	\subsection{Proof of Proposition~\ref{pro:main}} \label{app:main-proposition}
	
	\begin{proof}
		If $\varphi(z)\leq s$, the assumed local spectral gap gives the $(s,0)$-Poincar\'{e} inequality. Otherwise Popoviciu's inequality gives an $(s,1/4)$-Poincar\'{e} inequality, including when $\varphi(z)=\infty$.
		For part (i), Proposition~\ref{pro:conjugate-1}(i) gives the same spectral gap for $\bar K_0^*$; applying Lemma~\ref{lem:main}(i) gives the stated constants. For part (ii), positivity implies $\mathrm{Gap}(\bar K_2)\geq\mathrm{Gap}(\bar K_1)$, and Proposition~\ref{pro:conjugate-1}(i) and Lemma~\ref{lem:main}(ii) give the conclusion.
	\end{proof}
	
	\subsection{Proof of Corollary \ref{cor:uniform-gap}} \label{app:uniform-gap}
	
	\begin{proof}
		The first two bounds follow by applying Proposition~\ref{pro:main} with $s=1/\delta$ and $b_z=0$.
		We now prove the final bound.
		The argument is similar to the proof of Proposition 2.4 of \cite{CLY2026}, which establishes the bound for proximal samplers.
		See also \cite{qin2025spectral}.
		
		Let $(X,Y,Z)$ be distributed as $f_Z(\df z) \, f_{X \mid Z}(\df x \mid z) \, f_{Y \mid X}(\df y \mid x)$.
		Define the following DA Mtks:
		\[
		\begin{aligned}
			&K^{\ddagger}(y, \df y') = \int_{\X} f_{X \mid Y}(\df x \mid y) \, f_{Y \mid X}(\df y' \mid x), \\
			& \bar{K}_1^{\ddagger}(y, \df y') = \int_{\Z} f_{Z \mid Y}(\df z \mid y) \, f_{Y \mid Z}(\df y' \mid z), \\
			& K_z^{\ddagger}(y, \df y') = \int_{\X} f_{X \mid Y,Z}(\df x \mid y,z) \, f_{Y \mid X,Z}(\df y' \mid x,z),
		\end{aligned}
		\]
		where $f_{Y \mid X,Z}(\cdot \mid x,z) = f_{Y \mid X}(\cdot \mid x)$ almost surely.
		The Mtk $K^{\ddagger}$ defines a DA algorithm targeting $f_Y$, so does $\bar{K}_1^{\ddagger}$, while $K_z^{\ddagger}$ targets $f_{Y \mid Z}(\cdot \mid z)$.
		By Proposition \ref{pro:max-corr} in Appendix~\ref{app:da}, $\mathrm{Gap}(K^{\ddagger}) = \mathrm{Gap}(K)$, $\mathrm{Gap}(\bar{K}_1^{\ddagger}) = \mathrm{Gap}(\bar{K}_1)$, and $\mathrm{Gap}(K_z^{\ddagger}) = \mathrm{Gap}(K_z) \geq \delta$.
		
		For $h \in L_0^2(f_Y)$,
		\[
		\begin{aligned}
			\mathcal{E}_{K^{\ddagger}}(h) &= \frac{1}{2} \int f_Y(\df y) \, f_{Z \mid Y}(\df z \mid y) \, f_{X \mid Y,Z}(\df x \mid y,z) \, f_{Y \mid X}(\df y' \mid x) \, [h(y) - h(y')]^2 \\
			&= \frac{1}{2} \int_{\Y} f_Y(\df y) \int_{\Z} f_{Z \mid Y}(\df z \mid y) \int_{\Y} K_z^{\ddagger}(y, \df y') [h_z(y) - h_z(y')]^2 \\
			&= \int_{\Z} f_Z(\df z) \, \mathcal{E}_{K_z^{\ddagger}}(h_z) \\
			&\geq \delta \int_{\Z} f_Z(\df z) \|h_z\|_{f_{Y \mid Z}(\cdot \mid z)}^2,
		\end{aligned}
		\]
		where $h_z(\cdot) = h(\cdot) - \int_{\Y} h(y) \, f_{Y \mid Z}(\df y \mid z)$.
		The final line equals
		\[
		\delta \left\{ \|h\|_{f_Y}^2 - \int_{\Z} f_Z(\df z) \, \left[ \int_{\Y} h(y) \, f_{Y \mid Z}(\df y \mid z) \right]^2 \right\} = \delta \mathcal{E}_{\bar{K}_1^{{\ddagger}}}(h).
		\]
		Thus,
		\[
		\mathrm{Gap}(K) = \mathrm{Gap}(K^{\ddagger}) \geq \delta \mathrm{Gap}(\bar{K}_1^{\ddagger}) = \delta \mathrm{Gap}(\bar{K}_1).
		\]
	\end{proof}
	
	\subsection{Proof of Proposition~\ref{pro:converse}} \label{app:converse}
	
	\begin{proof}
		Assume that $\bar{K}_0$ satisfies an $(s,b)$-Poincar\'{e} inequality, where $s \geq 0$ and $b \geq 0$.
		For $h\in L_0^2(f_Z)$ and $A_0h(x)=\mathbb E[h(Z)\mid X=x]$, one has
		\[
		\mathcal E_{\bar K_1}(h)-\mathcal E_{\bar K_0}(h)
		=\langle A_0h,(I-K)A_0h\rangle_{f_X}\geq0.
		\]
		Thus $\bar K_1$ satisfies the same inequality. This argument uses reversibility and contraction of $K$ and does not require positivity.
		
		Assume instead that $K$ satisfies an $(s,b)$-Poincar\'{e} inequality, where $s \geq 1$ and $b \geq 0$.
		Let $g \in L_0^2(f_Z)$ be arbitrary.
		Then
		\[
		\begin{aligned}
			\int_{\Z^2} f_Z(\df z) \, \bar{K}_1(z, \df z') \, g(z) \, g(z') &= \int f_Z(\df z) f_{X \mid Z}(\df x \mid z) \, K(x, \df x') \, f_{Z \mid X}(\df z' \mid x') \, g(z) \, g(z') \\
			&= \int_{\X^2} f_X(\df x) \, K(x, \df x') \, l_g(x) \, l_g(x'),
		\end{aligned}
		\]
		where $l_g(\cdot) = \int_{\Z} f_{Z \mid X}(\df z \mid \cdot) \, g(z) \in L_0^2(f_X)$.
		Then, by the $(s,b)$-Poincar\'{e} inequality for~$K$ and Jensen's inequality,
		\[
		\begin{aligned}
			(1-s) \|g\|_{f_Z}^2 + s \langle g, \bar{K}_1 g \rangle_{f_Z} &= (1-s) \int_{\X \times \Z} f_X(\df x) \, f_{Z \mid X}(\df z \mid x) \, g(z)^2 + s \langle l_g, Kl_g \rangle_{f_X} \\
			&\leq (1-s) \|l_g\|_{f_X}^2 + s \langle l_g, K l_g \rangle_{f_X} \\
			&= \|l_g\|_{f_X}^2 - s \mathcal{E}_K(l_g) \\
			&\leq b\|l_g\|_{\mathrm{osc}}^2.
		\end{aligned}
		\]
		It is easy to see that $\|l_g\|_{\mathrm{osc}} \leq \|g\|_{\mathrm{osc}}$.
		It then follows from the most recent display that $\bar{K}_1$ satisfies an $(s,b)$-Poincar\'{e} inequality.
	\end{proof}

	\subsection{Proof of Proposition~\ref{pro:Kbar0-location}} \label{app:Kbar0-location}
	\begin{proof}
		By Proposition~\ref{pro:peskun-DA-0}, it suffices to show that, whenever $z \neq z'$,
		\[
		\exp\left[- \frac{L}{4} \|m(z) - m(z')\|^2 \right] \leq \left[ \int_{\mathbb{R}^d} \sqrt{ \exp\left[- \xi(x - m(z)) - \xi(x - m(z'))\right] } \, \df x \right]^2.
		\]
		To establish this, fix $(z,z') \in \Z^2$ such that $z \neq z'$, and let $\bar{m} = [m(z) + m(z')]/2$.
		Apply \eqref{ine:location-smooth} and obtain
		\[
		-\frac{1}{2} \xi(x - m(z)) - \frac{1}{2} \xi(x - m(z')) + \xi(x - \bar{m}) \geq - \frac{L\|m(z) - m(z')\|^2}{8}.
		\]
		Then, by the fact that $x \mapsto \exp[-\xi(x)]$ is a probability density function,
		\[
		\begin{aligned}
			&\int_{\mathbb{R}^d} \sqrt{ \exp\left[- \xi(x - m(z)) - \xi(x - m(z'))\right] } \, \df x \\
			&\geq \exp\left[-\frac{L\|m(z) - m(z')\|^2}{8}\right] \int_{\mathbb{R}^d} \exp[-\xi(x-\bar{m})] \, \df x \\
			&= \exp\left[-\frac{L\|m(z) - m(z')\|^2}{8}\right] 
		\end{aligned}
		\]
		as desired.
	\end{proof}
	
	\subsection{Proof of Proposition~\ref{pro:Kbar0-scale}} \label{app:Kbar0-scale}

	We start with a lemma regarding the Kullback-Leibler (KL) divergence between scaled versions of a centered, twice-differentiable, strongly log-concave, and log-smooth distribution.
	\begin{lemma} \label{lem:KL-scale}
		In the context of Section~\ref{sssec:Kbar0-scale}, when $0 < z < z'$, it holds that
		\[
		\mathrm{KL}(f_{X \mid Z}(\cdot \mid z), f_{X \mid Z}(\cdot \mid z')) \leq \frac{d}{8} \left( 1 + \frac{L}{a} \right) \left( \log z - \log z' \right)^2,
		\]
		where $\mathrm{KL}(\cdot, \cdot)$ is the KL divergence.
	\end{lemma}
	
	\begin{proof}
		Assume that $0 < z < z'$.
		Then, since the eigenvalues of $\mathrm{Hess}(\xi)$ are upper bounded by~$L$,
		\[
		\begin{aligned}
			&\mathrm{KL}(f_{X \mid Z}(\cdot \mid z), f_{X \mid Z}(\cdot \mid z')) \\
			=& \int_{\mathbb{R}^d} f_{X \mid Z}(\df x \mid z) \left[ -\xi\left( \frac{x}{\sqrt{z}} \right) + \xi\left( \frac{x}{\sqrt{z'}} \right) - \frac{d}{2} \left( \log z - \log z' \right)  \right] \\
			=& -\frac{d}{2} \, \log \frac{z}{z'} + \int_{\mathbb{R}^d} e^{-\xi(x)} \left[ -\xi(x) + \xi\left( \sqrt{\frac{z}{z'}} x \right) \right] \, \df x \\
			\leq& -\frac{d}{2} \, \log \frac{z}{z'} + \int_{\mathbb{R}^d} e^{-\xi(x)} \left[ \left( \sqrt{\frac{z}{z'}} - 1 \right) \nabla \xi(x)^{\top} x + \frac{L}{2} \left( \sqrt{\frac{z}{z'}} - 1 \right)^2 \|x\|^2 \right] \, \df x.
		\end{aligned}
		\]
		By Stein's identity \citep[][Lemma 2.2]{liu2016kernelized},
		the integral of $e^{-\xi(x)} \nabla \xi(x)^{\top} x$ over $\mathbb{R}^d$ is~$d$.
		On the other hand, since the eigenvalues of $\mathrm{Hess}(\xi)$ are lower bounded by~$a$, by the Brascamp-Lieb inequality \citep{brascamp1976extensions}, the integral of $e^{-\xi(x)} \, \|x\|^2$ over $\mathbb{R}^d$ is upper bounded by $d/a$.
		It follows that
		\[
		\mathrm{KL}(f_{X \mid Z}(\cdot \mid z), f_{X \mid Z}(\cdot \mid z')) \leq -\frac{d}{2} \, \log \frac{z}{z'} + d\left[ \exp\left( \frac{1}{2} \log \frac{z}{z'} \right) - 1 \right] + \frac{Ld}{2a} \left[ \exp\left( \frac{1}{2} \log \frac{z}{z'} \right) - 1 \right]^2.
		\]
		By a Taylor expansion of the exponential function, we have
		\[
		1 + \frac{1}{2} \, \log \frac{z}{z'} \leq \exp\left( \frac{1}{2} \log \frac{z}{z'} \right) \leq 1 \wedge \left[ 1 + \frac{1}{2} \, \log \frac{z}{z'} + \frac{1}{8} \left( \log \frac{z}{z'} \right)^2 \right].
		\]
		Thus,
		\[
		\begin{aligned}
			\mathrm{KL}(f_{X \mid Z}(\cdot \mid z), f_{X \mid Z}(\cdot \mid z')) \leq \frac{d}{8} \left( 1 + \frac{L}{a} \right) \left( \log \frac{z}{z'} \right)^2.
		\end{aligned}
		\]
	\end{proof}
	
	We now prove Proposition~\ref{pro:Kbar0-scale}.
	\begin{proof}
		Our proof relies on techniques involving isoperimetric inequalities and the conductance of a reversible Markov chain.
		
		If $f_{\log Z}$ is a point mass, then $\mathrm{Gap}(\bar K_0)=1$ by our convention, and the claimed bound holds.
		The bound is also trivial if $\mathrm{Is}(f_{\log Z})=0$, so assume that $f_{\log Z}$ is nondegenerate and $\mathrm{Is}(f_{\log Z})>0$.
		Define the metric between sets $\mathsf{d}_{\log}(\mathsf{S}_1,\mathsf{S}_2)=\inf_{z\in\mathsf{S}_1,\,z'\in\mathsf{S}_2} |\log z-\log z'|$ for measurable $\mathsf{S}_1, \mathsf{S}_2 \subset (0,\infty)$.
		
		Let $\nu(\df u)=\tfrac12 e^{-|u|}\,\df u$ be the standard Laplace distribution.
		By Theorem~1.3 of \cite{bobkov1997isoperimetric}, there exists a nondecreasing map $T:\mathbb{R}\to\mathbb{R}$ such that $\nu \circ T^{-1}=f_{\log Z}$ and $\|T\|_{\mathrm{Lip}}=1/\mathrm{Is}(f_{\log Z})$.
		The distribution $\nu$ has a positive density and isoperimetric profile $I_\nu(v)=\min\{v,1-v\}$, so Lemma~14 of \cite{andrieu2024explicit} gives a three set isoperimetric inequality, as we describe below.
		For a measurable partition $(0,\infty)=\mathsf{S}_1\sqcup\mathsf{S}_2\sqcup\mathsf{S}_3$ with $f_Z(\mathsf{S}_1),f_Z(\mathsf{S}_2)>0$,  take $\log \mathsf{S}_i = \{\log z: \, z \in \mathsf{S}_i\}$.
		Then
		\[
		\nu \circ T^{-1} (\log \mathsf{S}_3) \geq \mathrm{dist}(T^{-1}(\log \mathsf{S}_1), T^{-1}(\log \mathsf{S}_2)) \min \{\nu \circ T^{-1}(\log \mathsf{S}_1), \nu \circ T^{-1}(\log \mathsf{S}_2) \},
		\]
		where $\mathrm{dist}$ is the Euclidean distance between two sets.
		By the Lipschitz constant of~$T$,
		\[
		\operatorname{dist}(T^{-1}(\log \mathsf{S}_1),T^{-1}(\log \mathsf{S}_2))\geq \mathrm{Is}(f_{\log Z})\,\mathsf{d}_{\log}(\mathsf{S}_1,\mathsf{S}_2).
		\]
		This gives the following three-set isoperimetry:
		\begin{equation}\nonumber
		f_Z(\mathsf{S}_3)\geq \mathrm{Is}(f_{\log Z}) \,\mathsf{d}_{\log}(\mathsf{S}_1,\mathsf{S}_2)\min\{f_Z(\mathsf{S}_1),f_Z(\mathsf{S}_2)\}.
		\end{equation}

		On the other hand, $\bar{K}_0$ satisfies a so-called ``close coupling" condition.
		Let $z, z'$ be positive numbers.
		Then
		\begin{equation} \label{ine:TV-Kbar0-scale}
			\begin{aligned}
				\mathrm{TV}(\bar{K}_0(z, \cdot), \bar{K}_0(z',\cdot)) &= \sup_{g: \, (0,\infty) \to [0,1]} \left| \int_{\mathbb{R}^d} [f_{X \mid Z}(\df x \mid z) - f_{X \mid Z}(\df x \mid z')] \, \int_0^{\infty} f_{Z \mid X}(\df z'' \mid x) \, g(z'') \right| \\
				&\leq \sup_{h: \, \mathbb{R}^d \to [0,1]} \left| \int_{\mathbb{R}^d} [f_{X \mid Z}(\df x \mid z) - f_{X \mid Z}(\df x \mid z')] \, h(x) \right| \\
				&= \mathrm{TV}(f_{X \mid Z}(\cdot \mid z), f_{X \mid Z}(\cdot \mid z')).
			\end{aligned}
		\end{equation}
		By Lemma~\ref{lem:KL-scale} and Pinsker's inequality, if $z < z'$,
		\[
		\begin{aligned}
			\mathrm{TV}(f_{X \mid Z}(\cdot \mid z), f_{X \mid Z}(\cdot \mid z')) &\leq \frac{1}{4} \sqrt{d \left(1 + \frac{L}{a} \right)} |\log z - \log z'| .
		\end{aligned}
		\]
		Thus, by \eqref{ine:TV-Kbar0-scale} and symmetry, whenever $|\log z-\log z'|\leq 2(1+L/a)^{-1/2}d^{-1/2}$,
		\begin{equation} \nonumber
		\mathrm{TV}(\bar{K}_0(z, \cdot), \bar{K}_0(z',\cdot)) \leq \frac{1}{2}.
		\end{equation}
		This is referred to as a close coupling condition.
		
		A close coupling condition and three-set isoperimetry yields a conductance bound.
		This is a well-known fact; see, e.g., Lemma 15 of \cite{andrieu2024explicit} and Theorem 1.4.6 of \cite{qin2024convergence} for recent references.
		Applying either of the cited results gives
		\[
		\Phi:=\inf_{0<f_Z(\mathsf{A}) < 1}
		\frac{\int_{\mathsf{A}}\bar K_0(z,\mathsf{A}^c)\,f_Z(\df z)}{f_Z(\mathsf{A})f_Z(\mathsf{A}^c)}
		\geq \frac18\left[1\wedge\frac{\mathrm{Is}(f_{\log Z})}{\sqrt{(1+L/a)d}}\right].
		\]
		Since $\bar K_0$ is reversible, Cheeger's inequality \citep{lawler1988bounds} gives
		\[
		\mathrm{Gap}(\bar{K}_0) \geq \frac{\Phi^2}{2} \geq \frac{1}{128} \left[1 \wedge \frac{\mathrm{Is}(f_{\log Z})^2}{(1+L/a) d} \right].
		\]
	\end{proof}
	
	\subsection{Proof of Proposition~\ref{pro:Is-ig}} \label{app:Is-ig}

	\begin{proof}
		The density of $f_{\log Z}$ is
		\[
		Df_{\log Z}(u)=\frac{(\theta/2)^{\theta/2}}{\Gamma(\theta/2)}
		\exp\left(-\frac{\theta u}{2}-\frac{\theta e^{-u}}2\right),
		\qquad u\in\mathbb R.
		\]
		Its logarithm is concave. Consequently, for each $t\geq0$, the ratio
		$Df_{\log Z}(u-t)/Df_{\log Z}(u)$ is nondecreasing in $u$, whereas
		$Df_{\log Z}(u+t)/Df_{\log Z}(u)$ is nonincreasing.
		Integrating these ratios over $t\geq0$ gives
		\[
		\begin{aligned}
ote
			&\leq\frac{f_{\log Z}((-\infty,0])}{Df_{\log Z}(0)}
			\leq\frac1{Df_{\log Z}(0)}, &&u\leq0,\\
			\frac{f_{\log Z}((u,\infty))}{Df_{\log Z}(u)}
			&\leq\frac{f_{\log Z}((0,\infty))}{Df_{\log Z}(0)}
			\leq\frac1{Df_{\log Z}(0)}, &&u\geq0.
		\end{aligned}
		\]
		The definition of the isoperimetric constant therefore yields
		\begin{equation} \nonumber
			\mathrm{Is}(f_{\log Z})\geq Df_{\log Z}(0)
			=\frac{(\theta/2)^{\theta/2}\exp(-\theta/2)}{\Gamma(\theta/2)}.
		\end{equation}
		Using the gamma recurrence and the bound \citep[][Equation (3.9)]{artin1964gamma}
		\[
		\Gamma\left(\frac\theta2+1\right)
		\leq\sqrt{2\pi}
		\left(\frac{\theta+2}{2}\right)^{(\theta+1)/2}
		\exp\left(-\frac{\theta+2}{2}+\frac1{6(\theta+2)}\right),
		\]
		we obtain
		\[
		\begin{aligned}
			\mathrm{Is}(f_{\log Z})
			&\geq\frac{e}{2\sqrt\pi}\frac{\theta}{\sqrt{\theta+2}}
			\left(\frac\theta{\theta+2}\right)^{\theta/2}
			\exp\left(-\frac1{6(\theta+2)}\right)\\
			&\geq\frac{e^{-1/12}}{2\sqrt\pi}\frac{\theta}{\sqrt{\theta+2}}
			\geq\frac{\theta}{4\sqrt{\theta+2}},
		\end{aligned}
		\]
		where the second inequality uses $(1+2/\theta)^{\theta/2}\leq e$.
	\end{proof}

	\subsection{Proof of Corollary \ref{cor:mixing-MH-scale}} \label{app:mixing-MH-scale}
	
	\begin{proof}
		Let $s'' > 0$.
		In light of Proposition \ref{pro:Is-ig}, apply Proposition~\ref{pro:main}(i) with $s = s''/s'$, where
		\begin{equation} \nonumber
			s'=128\left[1\vee\frac{16(\theta+2)(1+L/a)d}{\theta^2}\right].
		\end{equation}
		Proposition~\ref{pro:Kz-MH-scale}, together with Lemma \ref{lem:ig-tail} stated right after this proof, 
		implies that $K$ satisfies an $(s'',b'')$-Poincar\'{e} inequality, with
		\[
		b'' =  \frac{s'}{4} \Pr \left( 10^4 \frac{Z}{\varsigma^2} e^{2\varsigma^2/Z} \frac{Ld}{a} > s \right) = \frac{s'}{4} \Pr \left( Z e^{2\varsigma^2/Z} > 10^{-4} \varsigma^2 \frac{a}{Ld} s \right) \leq \frac{s'}{4} c_{\theta} \left( 10^4 \varsigma^{-2} \frac{Ld}{a} \frac{s'}{s''} \right)^{\theta/2} ,
		\]
		where $c_{\theta} > 0$ satisfies $c_{\theta}^{2/\theta} = O(1)$ as $\theta \downarrow 0$.
		Letting $s''$ vary, we see that $K$ satisfies a weak Poincar\'{e} inequality with rate
		\[
		\beta(s'') = \frac{1}{4} c_{\theta} \left( 10^4 \varsigma^{-2} \frac{Ld}{a} \right)^{\theta/2} s'^{1+\theta/2} s''^{-\theta/2}.
		\]
		Corollary~\ref{cor:mixing-polynomial}, with $\ell=\theta/2$, implies that
		\[
		\begin{aligned}
			t(\epsilon; \mu, K) &\leq \left\lceil  c_{\theta}^{2/\theta} 10^4 \varsigma^{-2} (1+\theta/2)^{1+2/\theta} \frac{Ld}{a} s' \left( \frac{s'}{4} \frac{\|\mu - f_X \|_{f_X}^2}{\epsilon^2} \right)^{2/\theta} \right\rceil \\
			&= O \left( \frac{(L/a)(1+L/a)d^2}{\theta^2} \left[ \frac{512 (\theta+2) (1+L/a) d}{\theta^2}  \frac{\|\mu - f_X \|_{f_X}^2}{\epsilon^2} \right]^{2/\theta} \right).
		\end{aligned}
		\]
		Note that $(\theta+2)^{2/\theta} = 2^{2/\theta} (1+\theta/2)^{2/\theta} = O(2^{2/\theta})$ as $\theta \downarrow 0$.
		The desired mixing-time bound then follows.
	\end{proof}
	
	\begin{lemma} \label{lem:ig-tail}
		Suppose that $Z$ has the inverse gamma distribution with density 
		\begin{equation} \tag{\ref{eq:ig}}
			Df_Z(z) = \frac{(\theta/2)^{\theta/2}}{\Gamma(\theta/2)} z^{-\theta/2-1} \exp \left( - \frac{\theta}{2z} \right), \quad z > 0,
		\end{equation}
		where $\theta > 0$.
		Then, for $r > 0$,
		\begin{equation}\label{ine:ig-tail}
			\Pr(Z>r)
			\leq \frac{(\theta/2)^{\theta/2}}
			{\Gamma(1+\theta/2) \, r^{\theta/2}}.
		\end{equation}
		Moreover, for fixed $\varsigma \in (0,1/\sqrt{2})$ and $r > 0$,
		\begin{equation} \nonumber
			\begin{aligned}
				\Pr(Ze^{2\varsigma^2/Z}>r)
				\leq  c_{\theta} r^{-\theta/2},
			\end{aligned}
		\end{equation}
		where $c_{\theta} > 0$ satisfies $c_{\theta}^{2/\theta} = O(1)$ as $\theta \downarrow 0$.
	\end{lemma}
	
	\begin{proof}
		For every $r>0$, the inverse-gamma density gives
		\begin{equation} \nonumber
			\Pr(Z>r)
			=\frac1{\Gamma(\theta/2)}
			\int_0^{\theta/(2r)}u^{\theta/2-1}e^{-u}\,\df u
			\leq \frac1{\Gamma(\theta/2)}
			\int_0^{\theta/(2r)}u^{\theta/2-1} \,\df u
			= \frac{(\theta/2)^{\theta/2}}
			{\Gamma(1+\theta/2) \, r^{\theta/2}}.
		\end{equation}
		On $Z\geq1$, $Ze^{2\varsigma^2/Z}\leq e^{2\varsigma^2}Z$.
		On $Z<1$, 
		\[
		\mathrm E\!\left(Z^{\theta/2}e^{\theta\varsigma^2/Z}
		\ind_{\{Z<1\}}\right)
		=\frac{(\theta/2)^{\theta/2}}{\Gamma(\theta/2)}
		\int_1^\infty u^{-1}e^{-\theta(1/2-\varsigma^2)u}\,\df u.
		\]
		Combining~\eqref{ine:ig-tail} with Markov's inequality on the second event yields
		\begin{equation} \nonumber
			\begin{aligned}
				\Pr(Ze^{2\varsigma^2/Z}>r)
				\leq{}&\frac{(\theta/2)^{\theta/2}}
				{\Gamma(1+\theta/2)} \left[e^{\theta\varsigma^2}
				+\frac{\theta}{2}\int_1^\infty
				u^{-1}e^{-\theta(1/2-\varsigma^2)u}\,\df u\right] r^{-\theta/2} =: c_{\theta} r^{-\theta/2}.
			\end{aligned}
		\end{equation}
		This holds for every $\theta>0$.
		As $\theta\downarrow0$, the exponential integral in brackets is
		$-\log \theta + O(1)$ \citep{olver2026exponential}, while $- \log\Gamma(1+\theta/2)=O(\theta)$.
		Thus,
		\[
		c_{\theta}^{2/\theta} = \left\{
		\frac{(\theta/2)^{\theta/2}}{\Gamma(1+\theta/2)}
		\left[e^{\theta\varsigma^2}
		+\frac{\theta}{2}\int_1^\infty
		u^{-1}e^{-\theta(1/2-\varsigma^2)u}\,\df u\right]
		\right\}^{2/\theta}=O(1).
		\]
	\end{proof}
	
	\subsection{Proof of Proposition~\ref{pro:ordering-eta}} \label{app:ordering-eta}
	\begin{proof}
		For $0<\eta<\eta'$, couple the auxiliary variables as $Y_\eta=X+\sqrt\eta N$ and $Y_{\eta'}=Y_\eta+\sqrt{\eta'-\eta}N'$, where $N,N'$ are independent standard normal vectors, independent of $X$. Then $X,Y_\eta,Y_{\eta'}$ form a Markov chain. The conditional-expectation tower property and Jensen's inequality give, for every $g\in L_0^2(f_X)$,
		\[
		\mathbb E\{\mathbb E[g(X)\mid Y_{\eta'}]^2\}
		=\mathbb E\{\mathbb E[\mathbb E[g(X)\mid Y_\eta]\mid Y_{\eta'}]^2\}
		\leq \mathbb E\{\mathbb E[g(X)\mid Y_\eta]^2\}.
		\]
		Since $\mathcal E_{K(\eta)}(g)=\|g\|_{f_X}^2-\mathbb E\{\mathbb E[g(X)\mid Y_\eta]^2\}$, this proves $\mathcal E_{K(\eta')}(g)\geq\mathcal E_{K(\eta)}(g)$ and the desired ordering of Poincar\'{e} inequalities.
	\end{proof}
	\subsection{Proof of Proposition~\ref{pro:Kbar1-proximal-location}} \label{app:Kbar1-proximal-location}
	\begin{proof}
		Let $\phi_\eta$ be the density of $\mathrm N_d(0,\eta I_d)$. Fix $z,z'$ and set
		\[
		u=\frac{m(z')-m(z)}{1+L\eta},\qquad v=\frac{L\eta[m(z')-m(z)]}{1+L\eta}.
		\]
		Consider the product densities $e^{-\xi(r)}\phi_\eta(w)$ and $e^{-\xi(r-u)}\phi_\eta(w-v)$. Under the map $(r,w)\mapsto m(z)+r+w$, their pushforwards are $f_{Y\mid Z}(\cdot\mid z)$ and $f_{Y\mid Z}(\cdot\mid z')$, since $u+v=m(z')-m(z)$. Hellinger affinity increases under a measurable map, by conditional Cauchy--Schwarz. The product affinity factors, and the location estimate proved in Appendix~\ref{app:Kbar0-location}, together with the Gaussian affinity formula, gives
		\[
		\begin{aligned}
			\int\sqrt{Df_{Y\mid Z}(y\mid z)Df_{Y\mid Z}(y\mid z')}\,\df y
			&\geq \left[\int e^{-[\xi(r)+\xi(r-u)]/2}\,\df r\right]
			\left[\int\sqrt{\phi_\eta(w)\phi_\eta(w-v)}\,\df w\right]\\
			&\geq\exp\left[-\frac{L\|u\|^2}{8}-\frac{\|v\|^2}{8\eta}\right]
			=\exp\left[-\frac{L\|m(z')-m(z)\|^2}{8(1+L\eta)}\right].
		\end{aligned}
		\]
		Squaring and applying Proposition~\ref{pro:peskun-DA-1} proves the result.
	\end{proof}

	\subsection{Proof of Corollary \ref{cor:mixing-proximal-scale}} \label{app:mixing-proximal-scale}
	
	\begin{proof}
		Let $s'' > 0$.
		In light of Proposition \ref{pro:Is-ig}, apply Proposition~\ref{pro:main}(i) with $s = s''/s'$, where
		\begin{equation} \nonumber
			s'=128\left[1\vee\frac{16(\theta+2)(1+L/a)d}{\theta^2}\right].
		\end{equation}
		Proposition~\ref{pro:Kz-proximal-scale}, together with Lemma \ref{lem:ig-tail} in Appendix \ref{app:mixing-MH-scale},
		implies that $K$ satisfies an $(s'',b'')$-Poincar\'{e} inequality, with
		\[
		b'' =  \frac{1}{4} \wedge \frac{s'}{4} \Pr \left( \frac{Z+a\eta}{a\eta} > s \right) = \frac{1}{4} \wedge  \frac{s'}{4} \Pr \left( Z > a\eta \left( \frac{s''}{s'} - 1 \right) \right).
		\]
		By Lemma \ref{lem:ig-tail}, if $s'' > 2 s'$,
		\[
		b'' \leq \frac{s'}{4} \frac{(\theta/2)^{\theta/2}}{\Gamma(1+\theta/2)} \left[ a\eta \left( \frac{s''}{s'} - 1 \right) \right]^{-\theta/2} \leq \frac{s'}{4} \frac{(\theta/2)^{\theta/2}}{\Gamma(1+\theta/2)} \left( \frac{2s'}{a\eta s''} \right)^{\theta/2}.
		\]
		When $a\eta \leq 1$, for $\theta$ smaller than some universal constant, the right-hand-side is at least $1/4$ if $s'' \leq 2s'$.
		In such cases, we can claim $K$ satisfies a weak Poincar\'{e} inequality with rate 
		\[
		\beta(s'') = \frac{s'}{4} \frac{(\theta/2)^{\theta/2}}{\Gamma(1+\theta/2)} \left( \frac{2s'}{a\eta } \right)^{\theta/2} s''^{-\theta/2}.
		\]
		By Corollary~\ref{cor:mixing-polynomial},
		\[
		t(\epsilon; \mu, K) = O \left( \frac{(1+L/a)d}{a \eta \, \theta} \left[ \frac{512 (\theta+2) (1+L/a) d}{\theta^2 \Gamma(\theta/2 + 1)} \frac{\|\mu - f_X\|_{f_X}^2}{\epsilon^2} \right]^{2/\theta} \right).
		\]
		Note that, as $\theta \downarrow 0$, $(\theta+2)^{2/\theta} = O(2^{2/\theta})$, and $- \log \Gamma(\theta/2+1) = O(\theta)$.
		The desired result then follows.
	\end{proof}
	
	\subsection{Proof of Proposition~\ref{pro:location-semi}} \label{app:location-semi}
	\begin{proof}
		Write $U=-\log Df_X$, fix $x,h\in\mathbb R^d$, and let $\mathbb E_x$ denote expectation under $f_{Z\mid X}(\cdot\mid x)$. Then
		\[
		\nabla U(x)=\mathbb E_x[\nabla\xi(x-m(Z))].
		\]
		Set
		\[
		T_x=h^\top\nabla\xi(x-m(Z)),\qquad
		r_x=\xi(x+h-m(Z))-\xi(x-m(Z))-T_x.
		\]
		By~\eqref{ine:location-convex-smooth},
		\begin{equation} \label{ine:rx}
		\frac a2\|h\|^2\leq r_x\leq\frac L2\|h\|^2.
		\end{equation}
		Moreover, the $L$-Lipschitz property of $\nabla\xi$ gives
		\[
		\left|h^\top\nabla\xi(x-m(z))-h^\top\nabla\xi(x-m(z'))\right|
		\leq LR\|h\|.
		\]
		Thus, $T_x$ takes values in an interval of length at most $LR\|h\|$. Jensen's inequality and Hoeffding's lemma for moment generating functions imply
		\begin{equation} \label{ine:tx}
		0\leq\log\mathbb E_x \{\exp[-T_x + \mathbb E_x(T_x)]\}
		\leq\frac{L^2R^2}{8}\|h\|^2.
		\end{equation}
		The mixture density ratio satisfies
		\[
		\frac{Df_X(x+h)}{Df_X(x)}=\mathbb E_x [\exp(-T_x - r_x)],
		\]
		so
		\[
		U(x+h)-U(x)-\nabla U(x)^\top h
		=-\log\mathbb E_x \left\{\exp[-T_x+\mathbb E_x(T_x)-r_x] \right\}.
		\]
		Combining this with \eqref{ine:rx} and \eqref{ine:tx} yields
		\[
		\frac{a-L^2R^2/4}{2}\|h\|^2
		\leq U(x+h)-U(x)-\nabla U(x)^\top h
		\leq\frac L2\|h\|^2,
		\]
		which is~\eqref{ine:semi} with the claimed constants.
	\end{proof}
	
	\bibliographystyle{plainnat}
	\bibliography{qinbib}
\end{document}